\documentclass[11pt]{amsart}

\usepackage{amsmath, amsthm, amssymb, amscd, amsfonts}
\usepackage{fontsmpl}
\usepackage[all]{xy}

\newtheorem{lema}{Lemma}[section]
\newtheorem{theorem}[lema]{Theorem}
\newtheorem{cor}[lema]{Corollary}

\newtheorem{prop}[lema]{Proposition}

\theoremstyle{definition}
\newtheorem{definition}[lema]{Definition}

\theoremstyle{remark}
\newtheorem{obs}[lema]{Remark}

\newcommand\id{\operatorname{id}}

\newcommand\co{\operatorname{co}}

\newcommand\rk{\operatorname{rk}}

\newcommand\Alg{\operatorname{Alg}}
\newcommand\Aut{\operatorname{Aut}}
\newcommand\Int{\operatorname{Int}}
\newcommand\qAut{\operatorname{qAut}}
\newcommand\Map{\operatorname{Map}}

\newcommand\Ker{\operatorname{Ker}}
\newcommand\Img{\operatorname{Im}}

\newcommand\Emb{\operatorname{Emb}}

\newcommand{\eps}{\varepsilon}
\newcommand{\ot}{\otimes}
\newcommand{\com}{\Delta}

\newcommand\Res{\operatorname{Res}}
\newcommand\res{\operatorname{res}}
\newcommand{\g}{\mathfrak g}

\newcommand{\sop}{\operatorname{Supp}}

\newcommand{\D}{{\mathcal D}}

\newcommand{\Oc}{{\mathcal O}}

\newcommand\Fr{\operatorname{Fr}}
\newcommand\Lie{\operatorname{Lie}}

\def\ZZ {\mathbb{Z}}
\def\QQ{\mathbb{Q}}
\def\CC{\mathbb{C}}
\def\zh{\mathfrak{C}}
\def\GG{\mathbb{G}}

\def\JJ{\mathcal{J}}

\def\OO{\mathcal{O}}
\def\SS{\mathcal{S}}

\def\T{\mathbf{T}}
\def\te{\mathbb{T}}

\def\Z{\mathcal{HZ}}

\def\A{\mathcal{A}}

\def\D{\mathcal{D}}
\def\L{\mathcal{L}}
\def\J{\mathfrak{J}}
\def\t{\mathfrak{t}}

\def\lieg{\mathfrak{g}}
\def\lieh{\mathfrak{h}}

\def\liem{\mathfrak{m}}

\def\liesl{\mathfrak{sl}}

\def\Ga{\Gamma}
\def\Ge{G_{\epsilon}}
\def\Gli{\Ga(\lieg)}
\def\Glie{\Ga_{\epsilon}(\lieg)}

\def\Oe{\OO_{\epsilon}(G)}

\def\qe{{\bf u}_{\epsilon}(\lieg)}
\def\qsl{{\bf u}_{\epsilon}(\liesl_{2})}
\def\zh{\mathfrak{C}}
\def\liel{\mathfrak{l}}
\def\lgot{\mathfrak{l}}
\def\rgot{\mathfrak{r}}
\def\lo{\mathfrak{l}_{0}}

\def\Glil{\Ga_{\epsilon}(\liel)}
\def\Glir{\Ga_{\epsilon}(\rgot)}

\def\Oel{\OO_{\epsilon}(L)}

\def\qep{{\bf u}_{\epsilon'}(\lieg')}
\def\qsl{{\bf u}_{\epsilon}(\liesl_{2})}

\def\Ol{\Oc_{\epsilon}(L)}
\def\Or{\Oc_{\epsilon}(R)}

\def\ql{{\bf u}_{\epsilon}(\lgot)}
\def\qlo{{\bf u}_{\epsilon}(\lgot_{0})}
\def\qloe{{\bf u}_{\epsilon}(\lgot_{0}')}

\def\qr{{\bf u}_{\epsilon}(\rgot)}
\def\pf{\begin{proof}}
\def\epf{\end{proof}}
\theoremstyle{plain}

\def\tGa{\widetilde\Gamma}
\def\qlo{{\bf u}_{\epsilon}(\lgot_0)}

\begin{document}




\title[Extensions of finite quantum groups]
{Extensions of finite quantum groups by finite groups}

\author[n. andruskiewitsch and g. a. garc\'\i a]
{Nicol\' as Andruskiewitsch  \and Gast\' on Andr\' es Garc\'\i a}

\address{FaMAF-CIEM (CONICET), Universidad Nacional de C\'ordoba
Medina Allende s/n, Ciudad Universitaria, 5000 C\' ordoba, Rep\'
ublica Argentina.} \email{andrus@mate.uncor.edu,
ggarcia@mate.uncor.edu}

\thanks{Research of this paper was begun when G. A. G. was
visiting the Mathematisches Institut der Ludwig-Maximilians
Universit\" at M\"unchen under the support of the DAAD.  Results in
this paper are part of the Ph.D. thesis of G. A. G., written under
the advise of N. A. The work was partially supported by CONICET,
ANPCyT, Secyt (UNC) and Agencia Córdoba Ciencia}

\date{\today.}

\begin{abstract} We give a necessary and sufficient
condition for two Hopf algebras presented as central extensions to
be isomorphic, in a suitable setting. We then study the question of
isomorphism between the Hopf algebras constructed in \cite{AG} as
quantum subgroups of quantum groups at roots of 1. Finally, we apply
the first general result to show the existence of infinitely many
non-isomorphic Hopf algebras of the same dimension, presented as
extensions of finite quantum groups by finite groups.
\end{abstract}

\maketitle



\section{Introduction}
A major difficulty in the classification of finite-dimensional
Hopf algebras is the lack of enough examples, so that we are not
even able to state conjectures on the possible candidates to
exhaust different cases of the classification. Indeed, we are
aware at this time of the following examples of finite-dimensional
Hopf algebras: group algebras of finite groups; pointed Hopf
algebras with abelian group classified in \cite{AS-05} (these are
variations of the small quantum groups introduced by Lusztig
\cite{L1,L2}); other pointed Hopf algebras with abelian group
arising from the Nichols algebras discovered in \cite{G1, He}; a
few examples of pointed Hopf algebras with non-abelian group
\cite{MS, G2}; combinations of the preceding via standard
    operations: duals, twisting, Hopf subalgebras and quotients,
    extensions.

\smallbreak Let $G$ be a connected, simply connected, simple
complex algebraic group and let $\epsilon$ be a primitive
$\ell$-th root of 1, $\ell$ odd and $3 \nmid \ell$ if $G$ is of
type $G_{2}$. In \cite{AG}, we determined all Hopf algebra
quotients of the quantized coordinate algebra $\Oe$.
(Finite-dimensional Hopf algebra quotients of $\mathcal
O_{\epsilon}(SL_N)$ were previously obtained in \cite{Mu2}). A
byproduct of the main theorem of \cite{AG} is the discovery of
many new examples of Hopf algebras with finite dimension, or with
finite Gelfand-Kirillov dimension. The quotients in \cite{AG} are
parameterized by data $\D$ (see below for the precise definition)
and the corresponding quotient $A_{\D}$ fits into the following
commutative diagram with exact rows:

\begin{equation}\label{diag:const}
\xymatrix{1 \ar[r]^{} & \Oc(G) \ar[r]^{\iota} \ar[d]_{\res} & \Oe
\ar[r]^{\pi}\ar@/_2pc/@{.>}[dd]_(.7){q_{\D}}
\ar[d]^{\Res} & \qe^{*} \ar[r]^{}\ar[d]^{P} & 1\\
1 \ar[r]^{} &   \Oc(L)\ar[d]_{^{t}\sigma}  \ar[r]^{\iota_L} & \Ol
\ar[r]^{\pi_L}\ar[d]^{s} &
\ql^* \ar[r]^{}\ar[d] & 1\\
1 \ar[r]^{} &   \Oc(\Ga)  \ar[r]^{\hat{\iota}} &
A_{\D}\ar[r]^{\hat{\pi}}& H\ar[r]^{}& 1.}\end{equation}

The purpose of the present paper is to study when these examples
are really new (neither semisimple nor pointed, nor dual to
pointed), and when they are isomorphic as Hopf algebras. In
principle, it is hard to tell when two Hopf algebras presented by
extensions are isomorphic-- not as extensions but as `abstract'
Hopf algebras (even for groups there is no general answer). We
begin by studying in Section \ref{centralextensions} isomorphisms
between Hopf algebras of the form
\begin{equation}\label{diag:abstractsetting}
\xymatrix{ 1 \ar[r]^{} &   K  \ar[r]^{\iota} & A\ar[r]^{\pi}&
H\ar[r]^{}& 1,}\end{equation} where $K$ is the Hopf center of $A$,
and the Hopf center of $H$ is trivial. One of our main results--
Theorem \ref{isoytriple}--  gives a necessary and sufficient
condition for two Hopf algebras of this kind to be isomorphic, under
suitable hypothesis. Namely, we need $(i)$ $A$ noetherian and
$H$-Galois over $K$, and $(ii)$ any Hopf algebra automorphism of $H$
`lifts' to $A$. This setting is ample enough to include examples
arising from quantum group theory, and in particular from \cite{AG}.
Indeed, there is an algebraic group $\tGa$ such that $A_{\D}$ fits
into the following exact sequence
\begin{equation}\label{diag:centralext}
\xymatrix{ 1 \ar[r]^{} &   \Oc(\tGa)  \ar[r]^{\iota} &
A_{\D}\ar[r]^{\pi}& \qlo^*\ar[r]^{}& 1,}\end{equation} where
$\Oc(\tGa)$ is the Hopf center of $A_{\D}$, and the Hopf center of
$\qlo^*$ is trivial. As a first consequence, both $\tGa$ and $\qlo$
are invariants of the isomorphism class of $A_{\D}$, see Theorem
\ref{teo:invariantes}. However, the condition of lifting of
automorphisms remains an open question, except when $H = \ql^{*}$,
see Corollary \ref{propL}. Nevertheless, the Hopf center of
$\ql^{*}$ is trivial if and only if $\lgot = \lieg$. In this case,
we classify the quotients of $\Oe$ up to isomorphisms -- see Theorem
\ref{clasesiso}. Then using some results on cohomology of groups, we
prove that there are infinitely many non-isomorphic Hopf algebras of
the same dimension. They correspond to finite subgroup data and can
be constructed via a pushout. Using this fact, we are able to prove
that they form a family of non-semisimple, non-pointed Hopf algebras
with non-pointed duals. For $SL_{2}$ such an infinite family was
obtained by E. M\"uller \cite{Mu2}. Trying to understand this result
was one of our main motivations to study the problem of quantum
subgroups.

\subsection{Conventions}\label{conv}
Our references for the theory of Hopf algebras are \cite{Mo} and
\cite{Sw}, for Lie algebras \cite{Hu} and for quantum groups
\cite{J} and \cite{BG}. If $\Ga$ is a group, we denote by
$\widehat{\Ga}$ the character group. Let $k$ be a field. The
antipode of a Hopf algebra $H$ is denoted by $\mathcal S$.
\emph{All the Hopf algebras considered in this paper have
bijective antipode. See Remark \ref{skryabin}. } The Sweedler
notation is used for the comultiplication of $H$ but dropping the
summation symbol. The set of group-like elements of a coalgebra
$C$ is denoted by $G(C)$. We also denote by $C^{+} = \Ker \eps$
the augmentation ideal of $C$, where $\eps: C\to k$ is the counit
of $C$. Let $A\xrightarrow{\pi} H$ be a Hopf algebra map, then $
A^{\co H}= A^{\co \pi} =\{a \in A|\ (\id\otimes \pi)\Delta (a) =
a\otimes 1\}$ denotes the subalgebra of right coinvariants and $
^{\co H}A= \ ^{\co \pi}A =\{a \in A|\ (\pi\otimes \id)\Delta (a) =
1\otimes a\}$ denotes the subalgebra of left coinvariants.

\smallbreak A Hopf algebra $H$ is called {\it semisimple},
respectively {\it cosemisimple}, if it is semisimple as an algebra,
respectively if it is cosemisimple as a coalgebra. The sum of all
simple subcoalgebras is called the {\it{coradical}} of $H$ and it is
denoted by $H_{0}$. If all simple subcoalgebras of $H$ are
one-dimensional, then $H$ is called {\it pointed} and $H_{0} =
k[G(H)]$.

\smallbreak Let $H$ be a Hopf algebra,
$A$ a right $H$-comodule algebra
with structure map $ \delta: A\to A \ot H$,
$a \mapsto a_{(0)}\ot a_{(1)}$ and
$B = A^{\co H} $. The extension $B\subseteq A$
is called a Hopf Galois extension or $H$-Galois
if the canonical map $\beta: A\ot_{B} A\to A\ot H$,
$a\ot b \mapsto  ab_{(0)}\ot b_{(1)}$
is bijective. See \cite{SS} for more details on
$H$-Galois extensions.


\section{Central extensions of Hopf
algebras}\label{centralextensions}
\subsection{Preliminaries}\label{extension}
We recall some results on quotients and extensions of Hopf algebras.

\begin{definition}\cite{AD}\label{defsucex}
A sequence of Hopf algebras maps $1 \to B \xrightarrow{\iota} A
\xrightarrow{\pi} H \to 1$, where $1$ denotes the Hopf algebra $k$,
is {\it exact} if $\iota$ is injective, $\pi$ is surjective, $\Ker
\pi = AB^{+}$ and $B =\ ^{\co \pi} A$.
\end{definition}

\begin{obs}\label{obs:hgalois-ss}
Note that $A$ is a right $H$-Galois extension of $B$ by \cite{T},
see also \cite[3.1.1]{SS}.
\end{obs}

If the image of $B$ is central in $A$, then $A$ is called a {\it
central} extension of $B$. We say that $A$ is a {\it cleft
extension} of $B$ by $H$ if there is an $H$-colinear section
$\gamma$ of $\pi$ which is invertible with respect to the
convolution, see for example \cite[3.1.14]{A}. By \cite[Thm.
2.4]{Sch2}, a finite-dimensional Hopf algebra extension is always
cleft.

\smallbreak We shall use the following result.

\begin{prop}\cite[Prop. 2.10]{AG}\label{cociente1}
Let $A$ and $K$ be Hopf algebras, $B$ a central Hopf subalgebra of
$A$ such that $A$ is left or right faithfully flat over $B$ and
$p: B \to K$ a Hopf algebra epimorphism. Then $H = A/AB^{+}$ is a
Hopf algebra and $A$ fits into the exact sequence $1\to B
\xrightarrow{\iota} A \xrightarrow{\pi} H \to 1$. If we set $\JJ =
\Ker p \subseteq B$, then $(\JJ) = A\JJ $ is a Hopf ideal of $A$
and $A_{p} := A / (\JJ)$ is the pushout given by the following
diagram:

$$\xymatrix{B
\ar[r]^{\iota} \ar[d]_{p} & A \ar[d]^{q}\\
K\ar[r]_(.4){j} & A_{p} .}$$

\noindent Moreover, $K$ can be identified with a central Hopf
subalgebra of $A_{p}$ and $A_{p}$ fits into the exact sequence $
1\to K \to A_{p} \to H \to 1$. \qed
\end{prop}

\begin{obs}\label{decocaext} Let $A$, $B$ be as in Proposition
\ref{cociente1},
then the following diagram of central exact sequences is
commutative.
\begin{equation}\label{ext2} \xymatrix{1\ar[r]^{} &
B\ar[r]^{\iota} \ar[d]_{p} &
A\ar[r]^{\pi}\ar[d]^{q} & H\ar[r]^{}\ar@{=}[d]_{} & 1\\
1\ar[r]^{}& K\ar[r]^{j} &A_{p} \ar[r]^{\pi_{p}} & H \ar[r]^{}&
1.}\end{equation}
\end{obs}

\begin{obs}\label{apdimfinita}
If $\dim K$ and $\dim H$ is finite, then $\dim A_{p}$ is also
finite. Indeed, since the Galois map $\beta : A_{p}\otimes_{K}
A_{p} \to A_{p}\otimes H$ is bijective by Remark
\ref{obs:hgalois-ss} and $H$ is finite-dimensional, by \cite[Thm.
1.7]{KT}, $A_{p}$ is a finitely-generated projective $K$-module;
in particular $\dim A_{p}$ is finite.
\end{obs}

The following general lemma was kindly communicated to us by Akira
Masuoka.

\begin{lema}\label{lema:iso-galois-extensions}
Let $H$  be a bialgebra over an arbitrary commutative ring, and let
$A$,  $A'$  be right $H$-Galois extensions over a common algebra $B$
of $H$-coinvariants. Assume that $A'$  is right $B$-faithfully flat.
Then any $H$-comodule algebra map $\theta : A \to A'$ that is
identical on $B$  is an isomorphism. \end{lema}

\pf See \cite[Lemma 1.14]{AG}. \epf

\begin{obs}\label{obs:fivelema}
Masuoka's Lemma \ref{lema:iso-galois-extensions} implies the
following fact: let $A$ and $A'$ be Hopf algebra extensions of $B$
by $H$ and suppose that there is a Hopf algebra map $\theta: A\to
A'$ such that the following diagram commutes
$$\xymatrix{ 1 \ar[r]^{} & B \ar[r]^{}
\ar@{=}[d]_{} & A \ar[r]^{} \ar[d]^{\theta}
& H \ar[r]^{}\ar@{=}[d]^{} & 1\\
1 \ar[r]^{} &   B  \ar[r]^{} & A' \ar[r]^{} & H \ar[r]^{} & 1.}$$

\noindent If $A'$ is right $B$-faithfully flat, then $\theta$ must
be an isomorphism; cf. Remark \ref{obs:hgalois-ss}.
\end{obs}

\smallbreak The following proposition is due to E. M\" uller.

\begin{prop}\label{cociente2}\cite[3.4 (c)]{Mu2}
Let $1 \to B \xrightarrow{\iota} A \xrightarrow{\pi} H \to 1$ be
an exact sequence of Hopf algebras. Let $J$ be a Hopf ideal of $A$
of finite codimension and $\JJ = B\cap J$. Then $1\to B/ \JJ \to
A/ J \to  H / \pi(J) \to 1$ is exact. \qed
\end{prop}

\subsection{Isomorphisms}\label{isoobtext} Now we study some
properties of the Hopf algebras given by Proposition
\ref{cociente1}.

\begin{definition}\cite[2.2.3]{A}
The {\it Hopf center} of a Hopf algebra $A$ is the maximal central
Hopf subalgebra $\Z(A)$ of $A$. It always exists by \cite[2.2.2]{A}.
\end{definition}

\begin{prop}\label{isoinducidos}
For $i= 1,2$, let $1\to K_{i} \to A_{i} \to H_{i} \to 1$ be exact
sequences of Hopf algebras such that $K_{i}= \Z(A_{i})$. Suppose
that $\omega: A_{1}\to A_{2}$ is a Hopf algebra isomorphism. Then
there exist isomorphisms $\underline{\omega}: K_{1} \to K_{2}$ and
$\overline{\omega}: H_{1} \to H_{2}$ such that the following diagram
commutes
$$\xymatrix{ 1 \ar[r]^{} & K_{1} \ar[r]^{\iota_{1}}
\ar[d]_{\underline{\omega}} & A_{1} \ar[r]^{\pi_{1}} \ar[d]^{\omega}
& H_{1} \ar[r]^{}\ar[d]^{\overline{\omega}} & 1\\
1 \ar[r]^{} &   K_{2}  \ar[r]^{\iota_{2}} & A_{2} \ar[r]^{\pi_{2}} &
H_{2} \ar[r]^{} & 1.}$$
\end{prop}

\begin{proof}
Straightforward; see \cite[2.3.12]{G} for details.
\end{proof}

The following lemma and its corollaries will be needed later.

\begin{lema}\label{lema:centro}
Let $1\to B \to A \xrightarrow{\pi} H \to 1$ be a central exact
sequence of Hopf algebras such that $A$ is a noetherian. If $B
\subseteq C \subseteq \Z(A)$ is a Hopf subalgebra such that $\pi
(C) = \Z(H)$, then $C = \Z(A)$.
\end{lema}

\pf Let $B \subseteq D \subseteq \Z(A)$ be a Hopf subalgebra such
that $\pi (D) = \Z(H)$. By \cite[Thm. 3.3]{Sch1}, $A$ is
faithfully flat over $D$. Hence $D$ is a direct summand of $A$ as
$D$-module, see for example \cite[3.1.9]{SS}. Say $A = D\oplus M$.
Then $\Ker \pi|_{D} = DB^{+}$, since $\Ker \pi|_{D} = \Ker \pi
\cap D = AB^{+}\cap D= (D\oplus M)B^{+}\cap D = DB^{+}$. Besides
$B \subseteq D^{\co \pi|_{D}} \subseteq A^{\co \pi} = B$, which
implies that $D$ fits into the central exact sequence $ 1\to B \to
D \to \Z(H) \to 1$.

\smallbreak Moreover, the extension $B \subseteq D
\twoheadrightarrow \Z(H)$ is $\Z(H)$-Galois by Remark
\ref{obs:hgalois-ss}. Now taking $D=C$ and $D = \Z(A)$ we get the
following commutative diagram with exact rows which are
$\Z(H)$-Galois extensions of $B$
\begin{equation}\label{diag:five}
\xymatrix{1 \ar[r]^{} &   B \ar[r]^{}\ar@{=}[d] & C
\ar[r]^{}\ar@{^{(}->}[d]^{} & \Z(H)\ar[r]^{}\ar@{=}[d] & 1\\
1 \ar[r]^{} &  B  \ar[r]^{} & \Z(A) \ar[r]^{}& \Z(H)\ar[r]^{}& 1.}
\end{equation}

\noindent Hence, by Lemma \ref{lema:iso-galois-extensions} it
follows that $C = \Z(A)$. \epf

As immediate corollaries we get

\begin{cor}\label{cor:sincentro0} Assume the hypothesis of
Lemma \ref{lema:centro}. If $\Z(H) = k$, then $B = \Z(A)$.
\end{cor}

\pf It follows by taking $C= B$ in Lemma \ref{lema:centro}. \epf

\begin{cor}\cite[3.3.9]{A}\label{sincentro}
Let $1 \to K \xrightarrow{\iota}  A \xrightarrow{\pi} H\to1$ be a
central exact sequence of finite-dimensional Hopf algebras. If
$\Z(H) = k$, then $\Z(A) = K$. \qed
\end{cor}

We give now a sufficient condition for two Hopf algebras constructed
via the pushout to be isomorphic. Let $$\xymatrix{1 \ar[r]^{} & B
\ar[r]^{\iota} & A \ar[r]^{\pi} & H \ar[r]^{} & 1}$$

\noindent be a central exact sequence of Hopf algebras. Let $p_{1}:
B\to K_{1}$ and $p_{2}: B \to K_{2}$ be two Hopf algebras
epimorphisms. Then by Proposition \ref{cociente1}, we can build two
Hopf algebras $A_{1}:= A_{p_{1}}$ and $A_{2}:= A_{p_{2}}$, such that
$K_{i}$ is central in $A_{i}$ and by Remark \ref{decocaext}, both
fit into a commutative diagram:
\begin{equation}\label{diagcomm2} \xymatrix{1 \ar[r]^{} & B
\ar[r]^{\iota} \ar[d]_{p_{i}} & A \ar[r]^{\pi} \ar[d]^{q_{i}}
& H \ar[r]^{}\ar@{=}[d]_{} & 1\\
1 \ar[r]^{} & K_{i} \ar[r]^{j_{i}} & A_{i} \ar[r]^{\pi_{i}} &H
\ar[r]^{} & 1.}
\end{equation}

\begin{lema}\label{lema:isodeKs}
Let $f: K_{1} \to K_{2}$ be a Hopf algebra isomorphism such that $f
p_{1} = p_{2}$. Then the Hopf algebras $A_{1}$ and $A_{2}$ are
isomorphic.
\end{lema}

\begin{proof}
Straightforward; see \cite[2.3.14]{G} for details.
\end{proof}

We end this section with a theorem which gives under certain
assumptions  a characterization of the isomorphism classes of the
Hopf algebras obtained via the pushout construction. Note that Lemma
\ref{lema:isodeKs} gives a sufficient condition but under the
assumption that there is a Hopf algebra isomorphism $f: K_{1} \to
K_{2}$ such that $f p_{1} = p_{2}$. More generally, the theorem
below shows that it is enough to assume that the difference between
$f p_{1}$ and $p_{2}$ is given by a kind of dual of a 1-cocycle.

\smallbreak First, we need some definitions: we say that the central
$H$-extension $A$ of $B$ satisfies
\begin{enumerate}
\item[$(L)$] if every automorphism $f$ of $H$ can be
lifted to an automorphism $F$ of $A$ such that $\pi F = f \pi$, and
\item[$(Z)$] if $\Z(H) = k$.
\end{enumerate}

\bigbreak \emph{Assume for the rest of this section that $A$ is a
noetherian central extension of $B$. }

\bigbreak

\noindent Clearly, each Hopf algebra $A_{i}$ is noetherian, and
since $A_{i}$ is given by the pushout, one can see that
$K_{i}\subseteq A_{i}$ is an $H$-Galois extension. Indeed, let
$\beta: A\ot_{B} A\to A\ot H$, $a\ot b \mapsto ab_{(1)}\ot
\pi(b_{(2)})$ and $\beta_{i}: A_{i}\ot_{K_{i}} A_{i}\to A_{i}\ot H$,
$a'\ot b' \mapsto a'b_{(1)}'\ot \pi_{i}(b_{(2)}')$ be the canonical
maps. Since $A_{i} $ is a Hopf algebra extension of $K_{i}$ by $H$,
$\beta_{i}$ is surjective. Moreover, as $K_{i} = B / J_{i}$ and
$A_{i} = A/AJ_{i}$, $J_{i} = \Ker p_{i}$, the inverse of $\beta_{i}$
is given by $\beta_{i}^{-1}(q_{i}(a)\ot h) := \overline{(q_{i}\ot
q_{i})(\beta^{-1}(a\ot h)}$, for all $a \in A,\ h \in H$.

Let $f \in \Aut(H)$. If $(L)$ and $(Z)$ are satisfied, then by
Corollary \ref{cor:sincentro0}, $B = \Z(A)$ and $g= F|_{B}$ is an
automorphism of $B$. We denote by $\qAut(B)$ the subgroup of the
group $\Aut(B)$ of Hopf algebra automorphisms generated by these
elements.

\begin{theorem}\label{isoytriple}
Suppose that $(L)$ and $(Z)$ hold. Two Hopf algebras $A_{1}$ and
$A_{2}$ as in \eqref{diagcomm2} are isomorphic if and only if there
is a triple $(\underline{\omega},g, u)$ such that
\begin{enumerate}
\item[$(a)$] $\underline{\omega}: K_{1} \to K_{2}$ is an
isomorphism,
\item[$(b)$] $g \in \qAut(B)$,
\item[$(c)$] $u: A \to K_{2}$ is an algebra map and
\begin{eqnarray}
\underline{\omega}(p_{1}(b)) & = &
p_{2}(g(b_{(1)}))u(\iota(b_{(2)})),\label{conmutconu}\\
\com(u(a)) & = & u(a_{(2)})\otimes
q_{2}(F(\SS(a_{(1)})a_{(3)}))u(a_{(4)}),\label{comdeu}
\end{eqnarray}
\end{enumerate}
\noindent for all $b\in B$ and $a\in A$, where $F \in \Aut(A)$ is
induced by $\overline{\omega}$ with $F\iota = \iota g$.
\end{theorem}

\begin{obs} We are grateful to one of the anonymous referees for
the following interpretation of the conditions \eqref{conmutconu}
and \eqref{comdeu}. An isomorphism $\omega$ is identified with
\begin{equation}\label{referee}
A\ot_B K_1 \to A\ot_B K_2, \qquad a\ot x \mapsto F(a_{(1)})\ot
u(a_{(2)}) \underline{\omega}(x).
\end{equation}
Given a triple $(\underline{\omega},g, u)$, Condition
\eqref{conmutconu} (resp., \eqref{comdeu}) is equivalent to that
the (algebra) map \eqref{referee} is well-defined (resp., is a
coalgebra map).
\end{obs}

\begin{proof}
Let $\omega : A_{1} \to A_{2}$ be an isomorphism of Hopf algebras.
Since by assumption $\Z(H) = k$, from Corollary \ref{cor:sincentro0}
it follows that $\Z(A_{i}) = K_{i}$ for $1\leq i\leq 2$. Thus by
Proposition \ref{isoinducidos}, $\omega$ induces an isomorphism
$\underline{\omega}: K_{1} \to K_{2}$ and an automorphism
$\overline{\omega} \in \Aut(H)$. Then there exists an automorphism
$F$ of $A$ such that $\pi F = \overline{\omega} \pi$ and the map
given by $g= F|_{B}$ is an automorphism of $B$ such that $F \iota =
\iota g$. Define $u : A \to A_{2}$ to be the $k$-linear map given by
$$u(a) = q_{2}(F(\SS(a_{(1)})))\omega(q_{1}(a_{(2)}))\mbox{ for
all }a \in A,$$ that is, $u = q_{2}F\SS * \omega q_{1}$, the
convolution product between the maps $q_{2}F\SS$ and $\omega q_{1}$.
Since these maps are convolution invertible with inverses $q_{2}F$
and $\omega q_{1}\SS$ respectively,  $u$ is also convolution
invertible with inverse $\omega q_{1}\SS * q_{2}F$.

\bigbreak We claim that $u: A \to K_{2}$ is an algebra map and
satisfies \eqref{conmutconu} and \eqref{comdeu}. Indeed, it is
clear that $u(1) = 1$ and $\eps(u(a)) = \eps(a)$ for all $a \in
A$.

\bigbreak To prove that $\Img u \subseteq K_{2} =\ ^{\co\pi_{2}}
A_{2}$, let $a \in A$; then

$$\begin{array}{l} (\pi_{2}\otimes \id)\com(u(a)) = (\pi_{2}\otimes
\id) \com (q_{2}(F(\SS(a_{(1)})))\omega(q_{1}(a_{(2)})))
\\ [0.15cm]
= \pi_{2} (q_{2}(F(\SS(a_{(2)})))\omega(q_{1}(a_{(3)})))\otimes
q_{2}(F(\SS(a_{(1)})))\omega(q_{1}((a_{(4)})))
\\ [0.15cm]
= \pi_{2}
(q_{2}(F(\SS(a_{(2)})))\pi_{2}(\omega(q_{1}(a_{(3)})))\otimes
q_{2}(F(\SS(a_{(1)})))\omega(q_{1}((a_{(4)}))))
\\ [0.15cm]
=
\pi(F(\SS(a_{(2)})))\overline{\omega}\pi_{1}(q_{1}(a_{(3)}))\otimes
q_{2}(F(\SS(a_{(1)})))\omega(q_{1}((a_{(4)}))))
\\ [0.15cm]
=
\overline{\omega}\pi(\SS(a_{(2)}))\overline{\omega}\pi(a_{(3)})\otimes
q_{2}(F(\SS(a_{(1)})))\omega(q_{1}((a_{(4)}))))
\\ [0.15cm]
= \eps(a_{(2)})\otimes
q_{2}(F(\SS(a_{(1)})))\omega(q_{1}((a_{(3)}))))
\\ [0.15cm]
= 1\otimes q_{2}(F(\SS(a_{(1)})))\omega(q_{1}((a_{(2)})))) =
1\otimes u(a).
\end{array}$$

We prove now $u$ is an algebra map. Let $a,\ b \in A$, then

\begin{align*}
u(ab) & = q_{2}(F(\SS((ab)_{(1)})))\omega(q_{1}((ab)_{(2)}))\\
& = q_{2}(F(\SS(a_{(1)}b_{(1)})))\omega(q_{1}(a_{(2)}b_{(2)}))\\
& = q_{2}(F(\SS(b_{(1)})\SS(a_{(1)})))
\omega(q_{1}(a_{(2)}))\omega(q_{1}(b_{(2)}))\\
& = q_{2}(F(\SS(b_{(1)})))q_{2}(F(\SS(a_{(1)})))
\omega(q_{1}(a_{(2)}))\omega(q_{1}(b_{(2)}))\\
& = q_{2}(F(\SS(b_{(1)})))u(a)\omega(q_{1}(b_{(2)}))\\
& = u(a)(q_{2}(F(\SS(b_{(1)})))\omega(q_{1}(b_{(2)})))= u(a) u(b),
\end{align*}
\noindent since $u(a) \in K_{2}$ and $K_{2}$ is central in
$A_{2}$. Finally, let us prove that $u$ satisfies equations
\eqref{conmutconu} and \eqref{comdeu}: let $b \in B$, then
$u(\iota(b)) =
q_{2}(F\SS(\iota(b_{(1)})))\omega(q_{1}(\iota(b_{(2)})))$ and
therefore
\begin{align*}
\underline{\omega}(p_{1}(b)) & =
 \omega(q_{1}(\iota(b)))
= q_{2}(F(\iota(b_{(1)})))u(\iota(b_{(2)}))=
p_{2}(g(b_{(1)}))u(\iota(b_{(2)})).
\end{align*}

\noindent For the second equation, let $a \in A$, then
\begin{align*}
\com(u(a)) & = (q_{2}(F(\SS(a_{(1)})))\omega(q_{1}(a_{(2)})))_{(1)}
\otimes(q_{2}(F(\SS(a_{(1)})))\omega(q_{1}(a_{(2)})))_{(2)}\\
& = q_{2}(F(\SS(a_{(2)})))\omega(q_{1}(a_{(3)}))
\otimes q_{2}(F(\SS(a_{(1)})))\omega(q_{1}((a_{(4)})))\\
& = u(a_{(2)}) \otimes
q_{2}(F(\SS(a_{(1)})))\omega(q_{1}(a_{(3)}))\\
& = u(a_{(2)}) \otimes
q_{2}(F(\SS(a_{(1)})))q_{2}(F(a_{(3)}))u(a_{(4)})\\
& = u(a_{(2)}) \otimes
q_{2}(F(\SS(a_{(1)})a_{(3)}))u(a_{(4)}).\\
\end{align*}

Conversely, let $(\underline{\omega}, g, u)$ be a triple that
satisfies $(a)$, $(b)$ and $(c)$ and let $F \in \Aut(A)$,
$\overline{\omega} \in \Aut(H)$ such that $F|_{B} = g$ and $F \pi =
\overline{\omega}\pi$. Define $\varphi: A \to A_{2}$ to be the
$k$-linear map given by
$$\varphi(a) = q_{2}(F(a_{(1)}))u(a_{(2)})\mbox{ for all }a\in A.$$
\noindent As $K_{2}$ is central in $A_{2}$ and
 $u$ is an algebra map that
satisfies equation \eqref{comdeu}, it follows that $\varphi$ is a
Hopf algebra map. Moreover, by equation \eqref{conmutconu} we have
that
$$ \varphi(\iota(b)) = j_{2}(p_{2}(g(b_{(1)}))u(\iota(b_{(2)}))) =
j_{2}(\underline{\omega}(p_{1}(b)))\mbox{ for all }b \in B.$$

\noindent As $A_{1}$ is given by a pushout, there is a unique Hopf
algebra map $\omega: A_{1}\to A_{2}$ such that the following diagram
commutes:
$$\xymatrix{B
\ar[rr]^{\iota} \ar[d]_{p_{1}} & &
A \ar[d]_{q_{1}} \ar@/^/[ddr]^{\varphi}\\
 K_{1}\ar[rr]_{j_{1}} \ar@/_/[dr]_{\underline{\omega}}& & A_{1}
\ar@{-->}[dr]|-{\exists ! \omega}\\
& K_{2}\ar[rr]_{j_{2}} & & A_{2}.}$$

\noindent In particular, $ j_{2}\underline{\omega} = \omega j_{1}$
and $\overline{\omega}\pi_{1} = \pi_{2}\omega$, since for all $a \in
A$:
\begin{align*} \pi_{2}\omega (q_{1}(a)) & =
\pi_{2}\varphi (a) = \pi_{2} (q_{2}(F(a_{(1)}))u(a_{(2)})) =
\pi_{2}(q_{2}(F(a_{(1)}))) \pi_{2}(u(a_{(2)}))\\
& = \pi (F(a_{(1)})) \eps(a_{(2)}) = \overline{\omega}(\pi(a)) =
 \overline{\omega}(\pi_{1}(q_{1}(a)),
\end{align*}

\noindent and both exact sequences fit into a commutative diagram
\begin{equation}\label{diagcomm3} \xymatrix{1 \ar[r]^{} & K_{1}
\ar[r]^{j_{1}} \ar[d]_{\underline{\omega}} & A_{1} \ar[r]^{\pi_{1}}
\ar[d]^{\omega}
& H \ar[r]^{}\ar[d]^{\overline{\omega}} & 1\\
1 \ar[r]^{} & K_{2} \ar[r]^{j_{2}} & A_{2} \ar[r]^{\pi_{2}} &H
\ar[r]^{} & 1.}
\end{equation}

\noindent Since $\underline{\omega}$ and $\overline{\omega}$ are
Hopf algebra isomorphisms, the diagram above can be written as
\begin{equation}\label{diagcomm4} \xymatrix{1 \ar[r]^{} & K_{1}
\ar[r]^{j_{1}} \ar@{=}[d]_{} & A_{1} \ar[r]^{\pi_{1}}
\ar[d]^{\omega}
& H \ar[r]^{}\ar@{=}[d]^{} & 1\\
1 \ar[r]^{} & K_{1} \ar[r]^{j_{2}\underline{\omega}} & A_{2}
\ar[r]^{\overline{\omega}\pi_{2}} &H \ar[r]^{} & 1,}
\end{equation}

\noindent where the bottom exact sequence is exact. Since $A$ is a
noetherian $H$-Galois extension of $B$, then $A_{1}$ and $A_{2}$ are
noetherian, $K_{1} \subseteq A_{1}$ and $K_{2}\subseteq A_{2}$ are
$H$-Galois and by \cite[Thm. 3.3]{Sch1}, $A_{2}$ is faithfully flat
over $K_{2}$; in particular, $K_{1} \subseteq A_{2}$ is $H$-Galois
and $A_{2}$ is faithfully flat over $K_{1}$. Hence, by Remark
\ref{obs:fivelema}, $\omega$ is an isomorphism.
\end{proof}


\section{Quotients of the quantized coordinate algebra $\Oe$}\label{qca}
Let $G$ be a connected, simply connected complex simple Lie group
with Lie algebra $\lieg$, $\lieh\subseteq \g$ a fixed Cartan
subalgebra, $\Pi=\{\alpha_{1},\ldots,\alpha_{n} \}$ a basis of the
root system $\Phi = \Phi(\g,\lieh)$ of $\g$ with respect to
$\lieh$ and $n = \rk \g$. Let $\ell \geq 3$ be an odd integer,
relative prime to $3$ if $\lieg$ contains a $G_{2}$-component, and
let $\epsilon$ be a complex primitive $\ell$-th root of 1. By
$\Oe$ we denote the complex form of the quantized coordinate
algebra of $G$ at $\epsilon$ and by $\qe$ the Frobenius-Lusztig
kernel of $\lieg$ at $\epsilon$, see \cite{DL, BG} for
definitions. Let $\te :=\{K_{\alpha_{1}},\ldots ,K_{\alpha_{n}}\}
= G(\qe)$ be the ``finite torus'' of group-like elements of $\qe$
and for any subset $S$ of $\Pi$ define $\te_{S} : =
\{K_{\alpha_{i}}:\ \alpha_{i} \in S\}$.

\smallbreak We shall need the following facts about $\Oe$.

\begin{theorem}\label{OcentralenOe}
\begin{enumerate}
\item[$(a)$] \cite[Prop. 6.4]{DL}
$\Oe$ contains a central Hopf subalgebra isomorphic to the
coordinate algebra $\OO(G)$.
\item[$(b)$] \cite[III.7.11]{BG} $\Oe$ is a free
$\OO(G)$-module of rank $\ell^{\dim G}$.
\item[$(c)$] $\Oe$ fits into the following central exact sequence
\begin{equation}\label{extOOeqe}
 1\to \OO(G) \xrightarrow{\iota} \Oe \xrightarrow{\pi} \qe^{*} \to 1.
\end{equation}
\end{enumerate}
\end{theorem}

\pf $(c)$ follows immediately from $(a)$, $(b)$ and
\cite[III.7.10]{BG}. \epf

\begin{obs}\label{obs:referee} The extension \eqref{extOOeqe} is
now known to be cleft, see \cite[Th. 3.4.3]{SS}. Actually, $\pi$
admits a coalgebra section, see for example the proof of
\cite[Prop. 2.8 (c)]{AG}.
\end{obs}

A characterization of the Hopf algebra quotients of $\Oe$ is given
in \cite{AG}. We recall it briefly.

\begin{definition}\label{def:subgroupdatum}
We define a \emph{subgroup datum} for $G$ and $\epsilon$ as a
collection $\D= (I_{+}, I_{-}, N, \Ga, \sigma, \delta)$ where
\begin{itemize}
\item[$\bullet$] $I_{+}\subseteq \Pi$ and $I_{-} \subseteq -\Pi$.
Let $\Psi_{\pm} = \{\alpha \in \Phi:\ \sop \alpha \subseteq
I_{\pm}\}$, $\lgot_{\pm} = \sum_{\alpha \in \Psi_{\pm}} \g_{\alpha}$
and $\lgot = \lgot_{+} \oplus \lieh \oplus \lgot_{-}$; $\lgot$ is an
algebraic Lie subalgebra  of $\g$. Let $L$ be the connected Lie
subgroup of $G$ with $\Lie (L) = \lgot$. Let $I = I_{+} \cup -I_{-}$
and $I^{c} = \Pi - I$.

\medbreak \item[$\bullet$] $N$ is a subgroup of
$\widehat{\te_{I^{c}}}$.

\medbreak \item[$\bullet$] $\Ga$ is an algebraic group.

\medbreak \item[$\bullet$] $\sigma: \Ga \to L$ is an injective
homomorphism of algebraic groups.

\medbreak \item[$\bullet$] $\delta: N \to \widehat{\Ga}$ is a group
homomorphism.
\end{itemize}
\end{definition}

\noindent If $\Ga$ is finite, we call $\D$ a \emph{finite subgroup
datum}. The following theorem determines the Hopf algebra quotients
of $\Oe$.

\begin{theorem}\cite[Thm. 1]{AG}\label{teo:biyeccion}
There is a bijection between
\begin{enumerate}
\item[$(a)$] Hopf algebra quotients $q: \Oe \to A$. \item[$(b)$]
Subgroup data up to the equivalence relation defined in \cite[Def.
2.19]{AG}. \qed
\end{enumerate}
\end{theorem}

As an immediate consequence of Theorem \ref{teo:biyeccion}, there is
a bijection between Hopf algebra quotients $q: \Oe \to A$ such that
$\dim A < \infty$ and finite subgroup data up to equivalence.

\begin{obs}\label{skryabin} Given a subgroup datum $\D$, the
corresponding Hopf algebra quotient $A_\D$ fits into a central
exact sequence
\begin{equation*}
\xymatrix{ 1 \ar[r]^{} &   \Oc(\Ga)  \ar[r]^{\hat{\iota}} &
A_{\D}\ar[r]^{\hat{\pi}} & H \ar[r]^{} & 1,}\end{equation*} where
$H$ is the quotient of $\qe^{*}$ determined by the triple $(I_{+},
I_{-}, N)$ -- see \cite[Cor. 1.13]{AG}. Hence $A_\D$ satisfies the
hypothesis \textbf{(H)} in \cite[III.1.1, page 237]{BG}. By
\cite[Prop. III.1.1 part 1]{BG}, $A_\D$ is an affine noetherian PI
Hopf algebra. Therefore the antipode of $A_\D$ is bijective by
\cite[Cor. 2, page 623]{Sk}.
\end{obs}

\subsection{Some properties}\label{findimquo}
In this subsection we summarize some properties of the quotients
$A_{\D}$. Let $\D= (I_{+}, I_{-}, N, \Ga, \sigma, \delta)$ be a
subgroup datum. By \cite[Thm. 2.17]{AG}, $A_{\D}$ fits into the
commutative diagram

\begin{equation}\label{diag:facOel}
\xymatrix{1 \ar[r]^{} & \Oc(G)
\ar[r]^{\iota}\ar@/_2pc/@{.>}[dd]_(.7){^{t}\sigma} \ar[d]_{\res} &
\Oe \ar[r]^{\pi} \ar[d]^{\Res}\ar@/_2pc/@{.>}[dd]_(.7){q} & \qe^{*}
\ar[r]^{}\ar[d]_{p}
\ar@/^2pc/@{.>}[dd]^(.7){r} & 1\\
1 \ar[r]^{} &   \Oc(L)\ar[d]^{u}  \ar[r]^{\iota_L} & \Ol
\ar[r]^{\pi_L}\ar[d]^{w} &
\ql^* \ar[r]^{}\ar[d]_{v} & 1\\
1 \ar[r]^{} &   \Oc(\Ga)  \ar[r]^{\hat{\iota}} &
A_{\D}\ar[r]^{\hat{\pi}} & H \ar[r]^{} & 1,}\end{equation}

\noindent where $\lgot$ is an algebraic Lie subalgebra of $\lieg$,
$L$ the corresponding connected Lie subgroup of $G$ with $\Lie(L)
= \lgot$, $\ql^{*}$ is the quotient of $\qe^{*}$ determined by the
triple $(I_{+}, I_{-}, \te)$ and $\Oel$ the corresponding quantum
group. Moreover, by \cite[Lemma 2.14]{AG}, $H\simeq \ql^{*} /
(D^{z} - 1\vert\ z\in N)$.

\begin{lema}\label{lema:toro} The inclusion given by Theorem
\ref{OcentralenOe} $(a)$ determines a maximal torus $\T$ of $G$,
which is included in every connected subgroup $L$ corresponding to
the algebraic Lie subalgebra $\lgot \subseteq \lieg$ given by
Definition \ref{def:subgroupdatum}.
\end{lema}

\pf The first assertion is well-known, see \cite{DL} for details. We
sketch the proof to show that $\T$ is included in every $L$.

\smallbreak First, one defines an action of $\CC^{n}$ on $\Glie$:
let $\phi^{'}: \Gli \to \Glie$ be the canonical projection and
consider the primitive elements given by
$$H_{i} =
\phi^{'}\left(\frac{K_{\alpha_{i}}^{\ell}-1}{\ell(q_{i}^{\ell}-1)}
\right)
\in \Glie, \mbox{ for all }1\leq i\leq n.$$

\noindent By \cite[Def. 2.1]{AG}, these elements belong to the Hopf
subalgebra $\Glil$ of $\Glie$. Then for any $n$-tuple
$(p_{1},\ldots, p_{n}) \in \CC^{n}$ and for any finite-dimensional
$\Glil$-module $M$ the element $\exp(\sum_{i} p_{i}H_{i})$ defines
an operator which commutes with any $\Glil$-module homomorphism.
Hence, it defines a character on $\Oel$. Obviously, the elements
$\exp(\sum_{i} p_{i}H_{i})$ $\in \Oel^{*}$ form a group and the map
given by
\begin{align*}\phi: \CC^{n} \to
\Alg (\Oel, \CC), \quad (p_{1},\ldots, p_{n}) \mapsto
\exp\left(\sum_{i} p_{i}H_{i}\right),
\end{align*}
\noindent defines a group homomorphism whose kernel is the
subgroup $2\pi i\ell\ZZ^{n}$. Moreover, since by \cite[Thm.
10.8]{DL} the map $(\CC /2\pi i\ell\ZZ)^{n} \hookrightarrow \Alg
(\Oe, \CC)$ is an isomorphism and $\Alg (\Oel, \CC)
\hookrightarrow \Alg (\Oe, \CC)$, then the map $(\CC /2\pi
i\ell\ZZ)^{n} \hookrightarrow \Alg (\Oel, \CC)$ is also an
isomorphism.

Since the inclusion $\iota: \OO(G) \hookrightarrow \Oe$ given by
Theorem \ref{OcentralenOe} $(a)$ induces by restriction a group map
$ ^{t} \iota: \Alg (\Oe, \CC)\to \Alg (\OO(G), \CC)$, the
composition of this map with $\phi$ defines a homomorphism
$$\varphi: \CC^{n} \xrightarrow{\phi} \Alg (\Oe, \CC)\xrightarrow{^{t}
\iota}\Alg (\OO(G), \CC) = G,$$ \noindent whose image is contained
in $L = \Alg(\Oc(L), \CC)$ and its kernel is $2\pi i\ZZ^{n}$, by
\cite[Prop. 9.3 (c)]{DL}. The subgroup $\T\subseteq L \subseteq G$
given by the image of $\varphi$ is then a maximal torus of $G$. \epf

\begin{lema}\label{lema:dualnonpointed}
Denote by $j: q(\OO(G)) = \OO(\Ga) \to A_{\D}$ the inclusion map.
Then $j$ induces a group map $^{t}j: \Alg(A_{\D}, \CC) \to \Ga$
and $\Img (\sigma\circ\ ^{t}j) \subseteq \T\cap \sigma(\Ga)$.
\end{lema}

\begin{proof}
By \cite[Sec. 3]{AG}, the Hopf algebra $A_{\D}$ fits into the
following commutative diagram with exact rows
\begin{equation}\label{diag:dualpunteado}
\xymatrix{1 \ar[r]^{} & \OO(G) \ar[r]^{\iota} \ar[d]_{p} & \Oe
\ar[r]^{\pi} \ar[d]^{q}
& \qe^{*} \ar[r]^{}\ar[d]^{r} & 1\\
1 \ar[r]^{} &   \OO(\Ga)  \ar[r]^{j} & A_{\D} \ar[r]^{\bar{\pi}} & H
\ar[r]^{} & 1,}
\end{equation}

\noindent where $H$ is the finite-dimensional quotient of
$\qe^{*}$ determined by the triple $(I_{+}, I_{-}, N)$, see
\cite[Cor. 1.13]{AG}, \cite[Thm. 6.3]{MuI}. Then, the bottom exact
sequence induces an exact sequence of groups
\begin{equation}\label{extgrupos} 1 \to G(H^{*}) = \Alg(H, \CC)
\xrightarrow{^{t}\pi} \Alg(A_{\D}, \CC) \xrightarrow{^{t}j}
\Alg(\OO(\Ga), \CC) = \Ga,
\end{equation}
\noindent which fits into the commutative diagram of group maps
$$\xymatrix{1
\ar[r]^{} & G(\qe) \ar[r]^(.4){^{t}\pi} & \Alg(\Oe, \CC)
\ar[r]^(.45){^{t}\iota}
& \Alg(\OO(G),\CC) = G \\
1 \ar[r]^{} &   G(H^{*})  \ar[r]^{^{t}\pi}\ar[u]^{^{t}r} &
\Alg(A_{\D}, \CC) \ar[r]^{^{t}j}\ar[u]_{^{t}q} & \Ga
\ar[u]_{\sigma}.}$$

As $q$ is surjective, $^{t}q: \Alg(A_{\D}, \CC) \to \Alg(\Oe,
\CC)$ is injective. Since by Lemma \ref{lema:toro}, $\Alg(\Oe,
\CC) \simeq (\CC / 2\pi i \ell\ZZ)^{n}$ and the image of the
restriction map $^{t}\iota: \Alg(\Oe, \CC) \to G$ is the maximal
torus $\T$ of $G$, we have that the subgroup $(\sigma\circ\
^{t}j)(\Alg(A_{\D}, \CC))$ of $\sigma(\Ga)$ must be a subgroup of
$\T$.
\end{proof}

\smallbreak We resume some properties of the Hopf algebras
$A_{\D}$ in the following proposition. Part $(c)$ generalizes
\cite[Prop. 5.3]{Mu2}.

\begin{prop}\label{prop:propiedades}
Let $\D = (I_{+}, I_{-}, N, \Ga, \sigma, \delta)$ be a subgroup
datum.
\begin{enumerate}
  \item[$(a)$] If $A_{\D}$ is pointed, then $I_{+}\cap -I_{-} =
  \emptyset$ and $\Ga$ is a subgroup of the group of upper triangular matrices of some size.
  In particular, if $\Ga$ is finite, then it is abelian.
  \item[$(b)$] $A_{\D}$
  is semisimple if and only if $I_{+} \cup I_{-} = \emptyset$ and
  $\Ga$ is finite.
  \item[$(c)$] If $\dim A_{\D} < \infty$ and $A_{\D}^{*}$ is pointed, then
$\sigma(\Ga)\subseteq \T$. \item[$(d)$] If $A_{\D}$ is co-Frobenius
then $\Ga$ is reductive.
\end{enumerate}
\end{prop}

\pf $(a)$ If $A_{\D}$ is pointed, then by
\eqref{diag:dualpunteado} and \cite[Cor. 5.3.5]{Mo}, $H$ must be
also pointed. Thus $I_{+}\cap -I_{-} = \emptyset$, since otherwise
$H^{*}$ would contain a Hopf subalgebra isomorphic to $\qsl$ and
this would imply the existence of a surjective Hopf algebra map $H
\twoheadrightarrow \qsl^{*}$. This is impossible since $\qsl^{*}$
is non-pointed. Also, if $A_{\D}$ is pointed, then $\Oc(\Ga)$ is
also pointed, that is, any rational simple $\Ga$-module is
one-dimensional. Let $\rho: \Gamma \to GL(V)$ be a faithful
rational representation. Then $\rho(\Gamma)$ stabilizes a flag of
$V$, thus it is contained in a Borel subgroup of $GL(V)$.

\smallbreak $(b)$ If $\A_{D}$ is semisimple, then $H$ is also
semisimple and both are finite-dimensional. In such a case, $\Ga$
must be finite. Moreover, $H^{*}$ is also semisimple and this
implies by \cite[Thm. 6.3]{MuI} that $I_{+} \cup I_{-} = \emptyset$,
see also \cite[Cor. 1.13]{AG}. Conversely, if $\Ga$ is finite and
$I_{+} \cup I_{-} = \emptyset$, then $H$ is semisimple and whence
$A_{\D}$ is a central extension of two semisimple Hopf algebras,
implying as it is well-know that $A_{\D}$ is also semisimple.

\smallbreak $(c)$ If $\dim A_{\D}< \infty$, then the bottom exact
sequence of \eqref{diag:dualpunteado} induces the exact sequence
$1\to H^{*} \xrightarrow{^{t}\pi} A_{\D}^{*} \xrightarrow{^{t}j}
\CC[\Ga] \to 1$. Thus $^{t}j: A_{\D}^{*} \to \CC[\Ga]$ is surjective
and by \cite[Cor. 5.3.5]{Mo}, the image of the coradical of
$A_{\D}^{*}$ is the coradical of $\CC[\Ga]$, since $A_{D}^{*}$ is
pointed. Hence $^{t}j(G(A_{\D}^{*})) = \Ga$ and by Lemma
\ref{lema:dualnonpointed}, $\sigma(\Ga)$ must be a subgroup of $\T$.

\smallbreak $(d)$ A Hopf algebra $H$ is co-Frobenius if and only
if it is a right (or left) semiperfect coalgebra, see
\cite[5.3.2]{DNR}; and a subcoalgebra of a semiperfect coalgebra
is also semiperfect, see \cite[3.2.11]{DNR}. Hence, if $A_{\D}$ is
co-Frobenius, then $\Oc(\Ga)$ is co-Frobenius. Hence  $\Ga$ is
reductive by \cite{Su}. \epf

\begin{cor}\label{cor:pulenta}
Let $\D = (I_{+}, I_{-}, N, \Ga, \sigma, \delta)$ be a finite
subgroup datum such that $I_{+}\cap -I_{-} \neq \emptyset$ and
$\sigma(\Ga)\nsubseteq\T$. Then $A_{\D}$ is non-semisimple,
non-pointed and its dual is also non-pointed. \qed
\end{cor}

\subsection{Invariants}
Let $\qlo$ be the Hopf subalgebra of $\qe$ determined by the triple
$(I_{+}, I_{-}, \te_{I})$, see \cite[Cor. 1.13]{AG}. Then the
following properties concerning Hopf centers hold.

\begin{lema}\label{lema:qlo}
\begin{enumerate}
  \item[$(a)$] $\Z(\qlo^{*}) = \CC$.
  \item[$(b)$] $\qlo^{*} \simeq \ql^{*} /
(D^{z} - 1\vert\ z\in \widehat{\te_{I^{c}}})$.
  \item[$(c)$] $\Z(\ql^{*}) = \CC[\widehat{\te_{I^{c}}}]$.
  \item[$(d)$] $\Z(H) = \CC[\widehat{\te_{I^{c}}}/N]$.
\end{enumerate}
\end{lema}

\pf Since by \cite[Lemma A.2]{AS-app}, $\qlo$ is simple as a Hopf
algebra, it follows that $\Z(\qlo^{*}) = \CC$. By construction, it
is clear that $\qlo \subseteq H^{*}$ and $\qlo^{*} \simeq \ql^{*}
/ (D^{z} - 1\vert\ z\in \widehat{\te_{I^{c}}})$; this implies in
particular that $\qlo^{*}$ is a quotient of $H$.

\smallbreak By \cite[Lemma 2.14 (a)]{AG} we know that
$\CC[\widehat{\te_{I^{c}}}] \subseteq \Z(\ql^{*})$ and clearly
$\CC[\widehat{\te_{I^{c}}}]^{+}\ql^{*} = (D^{z} - 1\vert\ z\in
\widehat{\te_{I^{c}}})$. Thus, we have a central exact sequence of
finite-dimensional Hopf algebras
$$1\to \CC[\widehat{\te_{I^{c}}}] \to \ql^{*} \to \qlo^{*}\to 1,$$

\noindent with $\Z(\qlo^{*})=\CC$; in particular, $\ql^{*}$ is
$\qlo^{*}$-Galois by \cite{KT}. Then, by Corollary
\ref{cor:sincentro0}, $\Z(\ql^{*}) = \CC[\widehat{\te_{I^{c}}}]$
and part $(c)$ follows. The proof of $(d)$ follows directly from
the proof of $(c)$, since by \cite[Lemma 2.14 (b)]{AG}, $\qlo^{*}
\simeq H / (D^{z} - 1\vert\ z\in \widehat{\te_{I^{c}}}/N)$ and
$\CC[\widehat{\te_{I^{c}}}] \subseteq \Z(H)$. \epf

Given the algebraic group $\Ga$, we define now another algebraic
group $\widetilde{\Ga}$ which will help us to find invariants of
our construction.

\smallbreak From \cite[Prop. 2.6 $(b)$ and Lemma 2.17]{AG}, one
deduces that $\Oel$ contains a group of central group-like
elements isomorphic to $\widehat{\te_{I^{c}}}$. Thus,
$\Oc(L)\CC[\widehat{\te_{I^{c}}}]$ is a central Hopf subalgebra of
$\Oel$. Since it is commutative, there is an algebraic group
$\widetilde{L}$ such that $\Oc(L)\CC[\widehat{\te_{I^{c}}}] =
\Oc(\widetilde{L})$. Explicitly, $$\widetilde{L} =
\Alg(\Oc(L)\CC[\widehat{\te_{I^{c}}}],\CC) = L\times \te_{I^c}.$$

Analogously, by \cite[Lemma 2.17 and Thm. 2.18]{AG}, $A_{\D}$
contains a group of central group-like elements isomorphic to
$\widehat{\te_{I^{c}}}/N$ and therefore $\Oc(\widetilde{\Ga}):=
\Oc(\Ga)\CC[\widehat{\te_{I^{c}}}/N]$ is a central Hopf subalgebra
of $A_{\D}$, where $$\widetilde{\Ga} =
\Alg(\Oc(\Ga)\CC[\widehat{\te_{I^{c}}}/N],\CC) = \Ga\times
(\widehat{\te_{I^c}}/N)^{\widehat{}} = \Ga\times N^{\perp}$$

\noindent and $N^{\perp} = \{ x\in \te_{I_{c}}\vert\ \langle
z,x\rangle =0,\ \forall\ z\in N\}$.

\begin{lema}\label{lema:invariant}
\begin{enumerate}
  \item[$(a)$]  $\Z(\Oel) = \Oc(\widetilde{L})$
  and $\Oel$ fits into the central exact sequence
  $$1\to \Oc(\widetilde{L}) \to \Oel \to \qlo^{*} \to 1.$$
  \item[$(b)$] $\Z(A_{\D}) = \Oc(\widetilde{\Ga})$ and
  $A_{\D}$ fits into the central exact sequence
  $$1\to \Oc(\widetilde{\Ga}) \to A_{\D} \to \qlo^{*} \to 1.$$
  \item[$(c)$] $A_{\D}$ is the pushout given by the following
diagram:
$$\xymatrix{\Oc(\widetilde{L})
\ar[r]^{} \ar[d]_{} & \Oel \ar[d]^{}\\
\Oc(\widetilde{\Ga})\ar[r]_(.4){j} & A_{\D} .}$$
\end{enumerate}
\end{lema}

\pf $(a)$ Clearly $\Oc(\widetilde{L}) \subseteq \Z(\Oel)$. For the
other inclusion we apply Lemma \ref{lema:centro}: we know that
$\Oel$ fits into the central exact sequence
$$ 1\to \Oc(L) \to \Oel \to \ql^{*} \to 1,$$
which is also an $\ql^{*}$-Galois extension by \cite{KT}. Then the
lemma applies since $\Oel$ is noetherian and by Lemma
\ref{lema:qlo} $(c)$, $\pi_{L}(\Oc(L)\CC[\widehat{\te_{I^{c}}}])
=\CC[\widehat{\te_{I^{c}}}] = \Z(\ql^{*})$.

\smallbreak Recall that by Lemma \ref{lema:qlo} $(b)$, $\qlo^{*}
\simeq \ql^{*} / (D^{z} - 1\vert\ z\in \widehat{\te_{I^{c}}})$ and
denote by $v: \ql^{*}\to \qlo^{*}$ the surjective Hopf algebra map
given by the quotient. Then, we have the sequence of Hopf algebra
maps
\begin{equation}\label{suc:invariants}
    1\to \Oc(\widetilde{L}) \xrightarrow{j}
    \Oel \xrightarrow{v\pi_{L}} \qlo^{*} \to 1,
\end{equation}

\noindent with $j$ injective and $v\pi_{L}$ surjective. Moreover,
$\Ker v\pi_{L} = \Oel\Oc(\widetilde{L})^{+}$ and $\Oel^{\co
v\pi_{L}} = \Oc(L)$. Indeed, by definition it is clear that
$\Oc(\widetilde{L})^{+}\Oel \subseteq \Ker v\pi_{L}$, since
$\pi_{L}(\Oc(\widetilde{L})^{+}\Oel) =
\CC[\widehat{\te_{I^{c}}}]^{+} \ql^{*} = \Ker v$. Conversely, if
$a \in \Ker v\pi_{L}$, then $\pi_{L}(a)\in \Ker v =
\CC[\widehat{\te_{I^{c}}}]^{+}\ql^{*}$. Thus $a \in
\pi_{L}^{-1}(\CC[\widehat{\te_{I^{c}}}]^{+}\ql^{*})$ $=$ $\Ker
\pi_{L} + \CC[\widehat{\te_{I^{c}}}]^{+}\Oel =$ $\Oel\Oc(L)^{+}
+\CC[\widehat{\te_{I^{c}}}]^{+}\Oel \subseteq
\Oel\Oc(\widetilde{L})^{+}$, and this implies that $\Ker v\pi_{L}
\subseteq \Oel\Oc(\widetilde{L})^{+}$. Since $\Oel$ is a quotient
of $\Oe$, it is noetherian. Thus by \cite[Thm. 3.3]{Sch1}, $\Oel$
is faithfully flat over its central Hopf subalgebra
$\Oc(\widetilde{L})$. Since $\qlo^{*} \simeq \Oel /
[\Oc(\widetilde{L})^{+}\Oel]$, by \cite[Prop. 3.4.3]{Mo}, it
follows that $\Oel^{\co v\pi_{L}} = \Oc(\widetilde{L})$.

\smallbreak The proof of $(b)$ is analogous and we omit it. Just
take $H$ and $A_{\D}$ instead of $\ql^{*}$ and $\Oel$,
respectively.

\smallbreak $(c)$ Since $\Ga \subseteq L$, we have that
$\widetilde{\Ga} \subseteq \widetilde{L}$. Denote by $p:
\Oc(\widetilde{L}) \to \Oc(\widetilde{\Ga})$ the surjective Hopf
algebra map given by the transpose of this inclusion and let $K$ be
the Hopf algebra given by the pushout of $\Oc(\widetilde{L})
\hookrightarrow \Oel$ and $p$. Then we have the following
commutative diagram

$$\xymatrix{\Oc(\widetilde{L})
\ar@{^{(}->}[r]^{} \ar[d]_{p} &
\Oel \ar[d]^{} \ar@/^/[ddr]^{w}\\
\Oc(\widetilde{\Ga})\ar[r]_{} \ar@/_/[drr]_{j}& K
\\
& & A_{\D},}$$

\noindent where $w: \Oel \to A_{\D}$ is the map from
\eqref{diag:facOel}. Hence, there exists a Hopf algebra map
$\theta : K \to A_{\D}$ which makes the following commutative
diagram

$$\xymatrix{1 \ar[r]^{} & \Oc(\widetilde{\Ga}) \ar[r]^{} \ar@{=}[d]_{} &
K \ar[r]^{} \ar[d]^{\theta}
& \qlo^{*} \ar[r]^{}\ar@{=}[d]^{} & 1\\
1 \ar[r]^{} & \Oc(\widetilde{\Ga}) \ar[r]^{} & A_{\D} \ar[r]^{}
&\qlo^{*} \ar[r]^{} & 1.}
$$

\noindent Then, the proof of part $(c)$ follows by applying Remark
\ref{obs:fivelema}, since both extensions are $\qlo^{*}$-Galois by
\cite{KT}. \epf

Consequently, to any subgroup datum $\D$ we can attach an
algebraic group $\widetilde{\Ga}$ and an algebraic Lie algebra
$\lo$. The following theorem shows that these are invariants.

\begin{theorem}\label{teo:invariantes}
Let $\D$ and $\D'$ be subgroup data. If the Hopf algebras $A_{\D}$
and $A_{\D'}$ are isomorphic then $\widetilde{\Ga} \simeq
\widetilde{\Ga}'$ and $\lo \simeq \lo'$.
\end{theorem}

\pf Denote by $\theta: A_{\D} \to A_{\D'}$ the isomorphism. Then
$\theta(\Z(A_{\D})) = \Z(A_{\D'})$ and by Lemma
\ref{lema:invariant} $(a)$, the restriction of $\theta$ to
$\Oc(\widetilde{\Ga})$ defines an isomorphism of Hopf algebras
$\underline{\theta}: \Oc(\widetilde{\Ga}) \to
\Oc(\widetilde{\Ga}')$ and its transpose an isomorphism of
algebraic groups $\ ^{t}\theta: \Ga' \to \Ga$. Moreover, since by
Lemma \ref{lema:invariant} $(b)$, $\qlo^{*} = A_{\D}/
\Oc(\widetilde{\Ga})^{+}A_{\D}$ and $\qloe^{*} = A_{\D'}/
\Oc(\widetilde{\Ga}')^{+}A_{\D'}$, then $\theta$ also induces an
isomorphism $\overline{\theta}: \qlo^{*} \to \qloe^{*}$, which
implies that $\lo \simeq \lo'$. \epf



\section{An infinite family of Hopf algebras}
Now we focus in particular finite subgroup data to produce an
infinite family of non-isomorphic Hopf algebras.

\smallbreak Let $\Ga$ be an algebraic group and $\sigma: \Ga \to
G$ an injective homomorphism of algebraic groups. Denote by $p:
\OO(G) \to \OO(\Ga)$ the epimorphism of Hopf algebras given by the
transpose of $\sigma$. Then the exact sequence of Hopf algebras
\eqref{extOOeqe} gives rise by Proposition \ref{cociente1} to an
exact sequence
\begin{equation}\label{eqconstr1}
1\to \OO(\Ga) \xrightarrow{j} \Oe /(\JJ) \xrightarrow{\bar{\pi}}
\qe^{*} \to 1,
\end{equation}

\noindent where $\JJ = \Ker p$, $(\JJ) = \Oe\JJ$ and $\OO(\Ga)$ is
central in $\Oe /(\JJ)$. Thus, $A_{\sigma} := \Oe /(\JJ)$ is given
by a pushout. Moreover, since $A_{\sigma}$ is a quotient of $\Oe$
and $\qe^{*}$ is finite-dimensional, we have that $A_{\sigma}$ is a
noetherian $\qe^{*}$-Galois extension.

\smallbreak By Theorem \ref{teo:biyeccion}, this quotient of $\Oe$
corresponds to the subgroup datum $(\Pi, -\Pi, 1, \Ga, \sigma,
\eps)$, where $\eps: 1\to \widehat{\Ga}$ is the trivial group map.
If $\Ga$ is finite, then by Remark \ref{apdimfinita}, $\dim \Oe
/(\JJ) = \vert \Ga \vert \ell^{\dim \lieg}$ is also finite and
$A_{\sigma}$ corresponds to a finite subgroup datum. Furthermore, if
$\sigma(\Ga)\nsubseteq \T$, by Corollary \ref{cor:pulenta}
$A_{\sigma}$ is non-semisimple, non-pointed and its dual is also
non-pointed.

\subsection{The property $(L)$} Let $\qr\subseteq \qe$ be the
Frobenius-Lusztig kernel corresponding to a triple
$(I_{+},I_{-},\te_{\JJ})$ such that $\JJ\supseteq I= I_{+}\cup
-I_{-}$, and denote by $\rgot$ the corresponding Lie subalgebra of
$\g$; in particular, $\rgot \supseteq \lo$. Following \cite[Sec.
2]{AG}, one may  define the quantum groups $\Oc(R)$ and $\Or$; it
also holds that $\Or$ is a central $\qr^{*}$-extension of $\Oc(R)$.
We show that these kind of extensions satisfy the property $(L)$ of
Subsection \ref{isoobtext}. In particular, the results also hold for
the $\qe^{*}$-extension $\Oe$ of $\Oc(G)$.

\begin{lema}\label{indaut}
Every automorphism $f$ of $\qr$ induces an automorphism
of $R = \Alg(\Oc(R), \CC)$.
If moreover $\lieh \subseteq \rgot$, then $f$ leaves
invariant the torus $\T$.
\end{lema}

\begin{proof}
Let $\overline{F}$ be an automorphism of $\qr$. Since this Hopf
algebra is a quantum linear space of finite Cartan type, by
\cite[Thm. 7.2]{AS-05} and the proof of \cite[Thm. 5.9]{MuI},
there is a unique automorphism $F$ of $\Glir$ such that $F|_{\qr}
= \overline{F}$.

\smallbreak Consider now the quantum Frobenius map $\Fr: \Glie \to
U(\lieg)_{\QQ(\epsilon)}$ and denote by $\overline{\Fr}: \Glir \to
U(\rgot)_{\QQ(\epsilon)}$ its restriction, which is defined on the
generators of $\Glir$ by
\begin{align*}
\overline{\Fr}(E_{j}^{(m)}) =
\begin{cases}{\begin{matrix}e_{j}^{(m/\ell)} & \mbox{if }
\ell\vert m\\
0& \mbox{otherwise}\end{matrix}}\end{cases}, & &
\overline{\Fr}(F_{k}^{(m)}) =
\begin{cases}\begin{matrix}f_{k}^{(m/\ell)} & \mbox{if } \ell\vert
m\\
0& \mbox{otherwise}\end{matrix}\end{cases},\\
\overline{\Fr}(\left(\begin{smallmatrix} K_{\alpha_{i}};0\\
m\end{smallmatrix}\right))
= \begin{cases}\begin{matrix} \left(\begin{smallmatrix} h_{i};0\\
m\end{smallmatrix}\right)& \mbox{if } \ell\vert
m\\
0& \mbox{otherwise}\end{matrix}\end{cases}, & &
\overline{\Fr}(K_{\alpha_{i}}^{-1}) = 1,\qquad\qquad
\end{align*}

\noindent for all $j\in I_{+}$, $k\in I_{-}$ and $i\in \JJ$.
Then by \cite[Thm. 6.3]{DL}, $\Fr$ is a well-defined Hopf algebra
map and following \cite[Rmk. 2.5 $(b)$]{AG}, one sees that
the kernel of $\overline{\Fr}$
is the two-sided ideal generated by the set
$$\Big\{E_{j}^{(m)},F_{k}^{(m)}, \left(\begin{smallmatrix} K_{\alpha_{i}};0\\
m\end{smallmatrix}\right), K_{\alpha_{i}} -1, p_{\ell}(q) \vert\
i \in \JJ, j\in I_{+}, k\in I_{-}, \ell \nmid m\Big\}.$$

Using again the explicit description of the automorphism group of
$\qr$ given in \cite[Thm. 7.2]{AS-05}, it follows that
$F(\Ker\overline{\Fr}) = \Ker \overline{\Fr}$. Hence $F$
factorizes through a Hopf algebra automorphism
$\underline{F}:U(\rgot)_{\QQ(\epsilon)} \to
U(\rgot)_{\QQ(\epsilon)}$. Since by definition $\Oc(R)\subseteq
U(\rgot)^{\circ}$ -- see \cite[Sec. 2]{AG}, the transpose
$^{t}\underline{F}$ of $\underline{F}$ induces an automorphism of
$\OO(R)$, obtaining in this way an automorphism $f$ of $R =
\Alg(\Oc(R), \CC)$ which comes from an automorphism
$\underline{F}$ of $\qr$.

\smallbreak Suppose now that $\rgot \supseteq \lieh$, that is
$\te_{\JJ} = \te$. Then by Lemma \ref{lema:toro}, $\rgot$
corresponds to an algebraic Lie subalgebra of $\lieg$ which
contains the torus $\T$ and $\rgot = \Lie(R)$. We will show that
$f(\T) = \T$. Recall that by the proof of Lemma \ref{lema:toro},
$\T =\ ^{t}\iota (G_{\epsilon}) =\ ^{t}\iota_{R} (R_{\epsilon})$,
where $G_{\epsilon} = \Alg(\Oe, \CC)$ and $R_{\epsilon} =
\Alg(\Or, \CC)$. By \cite[4.1 and 6.1]{DL}, there is a perfect
pairing between $\Oe$ and $\Glie$, which clearly restricts to a
perfect pairing between $\Or$ and $\Glir$. Thus, the automorphism
$F$ induces an automorphism $^{t}F$ of $\Or$. Since $F$ factorizes
through $\underline{F}$ we have that $^{t}F \iota_{R} = \iota_{R}\
^{t}\underline{F}$ and hence,
\begin{align*}
f(\T) & = f(^{t}\iota_{R}(R_{\epsilon})) =\ ^{t}(\iota_{R}\
^{t}\underline{F})(R_{\epsilon}) =\ ^{t}(^{t}F
\iota_{R})(R_{\epsilon}) =\ ^{t}\iota_{R}(F(R_{\epsilon})) \\
& =\ ^{t}\iota_{R}(R_{\epsilon}) \subseteq \T .\end{align*}

\noindent Thus $f(\T) = \T$, because $f$ is an automorphism.
\end{proof}

\begin{cor}\label{propL}
The $\qr^{*}$-extension $\Or$ of $\OO(R)$ satisfies $(L)$.
\end{cor}

\begin{proof}
Since $\dim\qe<\infty$, every automorphism $\alpha$ of $\qr^{*}$
corresponds to an automorphism $\overline{F}$ of $\qr$. Thus, from
the proof of the lemma above, $\overline{F}$ induces an automorphism
$F$ of $\Glir$ such that $F ^{t}\pi_{R} =\ ^{t}\pi_{R}
\overline{F}$. Hence $^{t}F \in \Aut(\Or)$ and $\alpha \pi_{R} =
\pi_{R}\ ^{t}F$, which implies the claim.
\end{proof}

\begin{definition}
Denote by $\qAut(R)$ the subgroup of the group of rational
isomorphism $\Aut(R)$ of $R$ generated by isomorphisms of $R$ coming
from automorphisms of $\qr$.
\end{definition}

\subsection{The group $\qAut(G)$}

From now on, we assume that $R = G$. Let $B$ be the Borel subgroup
of $G$ that contains $\T$; it is determined by fixing the base
$\com$ of the root system $\Phi$ determined by $\T$. Let $D$ be
the subgroup of $\Aut(G)$ given by
$$D = \{ f\in \Aut(G)\vert\ f(\T) =  \T\mbox{ and } f(B) = B\}.$$

\noindent By \cite[Cor. 5.7]{MuI}, $D\subseteq\qAut(G)$. Since
$f(\T) = \T$, each $f \in D$ induces an automorphism $\hat{f}$ of
$\Phi$. Moreover, since $f(B) = B$, $f$ preserves $\com$ and whence
$\hat{f}$ belongs to the group of diagram automorphisms of $\Phi$.
If $\Int(G)$ denotes the subgroup of inner automorphisms of $G$,
then by \cite[Thm. 27.4]{Hu2}, $\Int(G)$ is normal in $\Aut(G)$ and
$\Aut(G) = \Int(G)\rtimes D$; in particular, $\Int(G)$ has finite
index in $\Aut(G)$. Since for all $t\in \T$, the inner automorphism
$\Int(t)$ of $G$ given by the conjugation fixes $\T$ and $B$, see
\cite[Lemma 24.1]{Hu2}, the image $\Int(\T)$ of $\T$ in $\Aut(G)$ is
a subgroup of $D$.

\smallbreak Denote by $\Int(N_{G}(\T))$ the subgroup of inner
automorphisms of $\Aut(G)$ given by the conjugation of elements in
the normalizer $N_{G}(\T)$ of $\T$ in $G$.

\begin{lema}\label{TactuaenqAut}
\begin{enumerate}
\item[$(a)$] $\qAut(G)$ is a subgroup of $\Int(N_{G}(\T))\rtimes
D$. \item[$(b)$] $\T$ acts on $\qAut(G)$ by left multiplication of
$\Int(\T)$. \item[$(c)$] The set $\qAut(G)/\T$ of orbits of the
preceding action is finite.
\end{enumerate}
\end{lema}

\begin{proof}
$(a)$ Let $f \in \qAut(G)$, then there exist $\alpha \in \Int(G)$,
$\beta \in D$ such that $f=\alpha \beta$. Since $f(\T) = \T$ and
$\beta(\T) = \T$, it follows that $\alpha(\T) = \T$. Hence $\alpha
= \Int(g)$, for some $g \in N_{G}(\T)$.

\smallbreak $(b)$ Since $\Int(\T)\subseteq D\subseteq\qAut(G)$,
the left multiplication by elements of $\Int(\T)$ defines an
action of $\T$ on $\qAut(G)$.

\smallbreak $(c)$ Since by $(a)$, $\qAut(G) \subseteq
\Int(N_{G}(\T)) \rtimes D$, it follows that
\begin{align*}
\vert \qAut(G) / \T \vert & \leq \vert [\Int(N_{G}(\T))\rtimes D]
/ \T \vert \leq \vert N_{G}(\T) / \T \vert \vert D\vert = \vert
W_{\T} \vert \vert D\vert,
\end{align*}

\noindent where $W_{\T} = N_{G}(\T) / C_{G}(\T)= N_{G}(\T) / \T$
is the Weyl group associated to $\T$. The claim follows since the
orders of $W_{\T}$ and $D$ are finite.
\end{proof}

\subsection{Group cohomology}
To describe the isomorphism classes of this type of extensions we
shall need some basic facts from  cohomology of groups. Let $G$,
$\Ga$ be two groups and suppose that there exists a right action $
\leftharpoonup$ of $\Ga$ on $G$ by group automorphisms. By $\Map
(\Ga, G)$ we denote the set of maps from $\Ga$ to $G$. For $n=0,1$
we define differential maps $\partial_{n}$ by
\begin{align*}
\partial_{0}: & \Map (1,G) \to \Map (\Ga, G), & \partial_{0}(g)(x) & =
(g\leftharpoonup x)g^{-1},\\
\partial_{1}:  & \Map (\Ga,G) \to \Map (\Ga\times \Ga, G),
& \partial_{1}(v)(x,y) & = (v(x)\leftharpoonup y)v(y)v(xy)^{-1},
\end{align*}
\noindent for all $x,\ y \in \Ga, g\in G$ and $v \in \Map (\Ga,G)$.
As usual, $\partial^{2} = 1$.

\begin{definition}
$(i)$ A map $u \in \Map(\Ga,G)$ is called a 1-{\it coboundary} if
$u \in \Img \partial_{0}$, that is, if there exists $g\in G$ such
that $u(x) = \partial_{0}(g)(x) = (g\leftharpoonup x)g^{-1}$ for
all $x\in \Ga$.

\noindent $(ii)$ A map $v \in \Map(\Ga,G)$ is called a 1-{\it
cocycle} if $\partial_{1}(v) = 1$, that is, if $v(xy) =
(v(x)\leftharpoonup y)v(y)$ for all $x,\ y \in \Ga$.
\end{definition}

Clearly, every 1-coboundary is a 1-cocycle. Denote by $Z^{1}(\Ga,
G)$ the subset of 1-cocycles in $\Map(\Ga,G)$. Then $G = \Map(1,G)$
acts on $Z^{1}(\Ga, G)$ via
\begin{equation} (g\cdot v)(x) =
(g\leftharpoonup x)v(x)g^{-1},
\end{equation}

\noindent for all $g\in G$, $v \in Z^{1}(\Ga, G)$ and $x\in \Ga$.
We say that two 1-cocycles $v$ and $u$ are {\it equivalent},
$v\sim u$, if there exists $g \in G$ such that $v = g\cdot u$.
Then we set $H^{1}(\Ga, G) := Z^{1}(\Ga, G)/G $. In particular,
$\bar{v} = \bar{1}$ in $H^{1}(\Ga, G)$ if and only if $v$ is a
1-coboundary.

\smallbreak Now we apply these notions in our setting. Let $G$ be a
connected, simply connected, semisimple complex Lie group as in
Section \ref{qca} and let $\sigma: \Ga \to G$ be an injective
homomorphism of algebraic groups. For any  $f \in \Aut(G)$ we define
an action of $\Ga$ on $G$, depending on $\sigma$ and $f$, via the
conjugation:
\begin{equation}\label{acciondegaeng} G\times \Ga
\xrightarrow{\leftharpoonup} G,\qquad g\leftharpoonup x  =
(f\sigma(x))^{-1} g (f\sigma(x)),
\end{equation}

\noindent for all $g\in G,\ x \in \Ga$. Hence,  $u \in \Map(\Ga,G)$
is a 1-coboundary if and only if there exists $g \in G$ such that
\begin{equation}\label{cocadena} u(x) =
\partial_{0}(g)(x) = (g\leftharpoonup x)g^{-1} = (f\sigma(x))^{-1}
g (f\sigma(x)) g^{-1},
\end{equation}

\noindent for all $x \in \Ga$, and a map $v\in \Map(\Ga, G)$ is a
1-cocycle if and only if
\begin{equation}\label{cociclo} v(xy) =
(v(x)\leftharpoonup y) v(y) = (f\sigma(y))^{-1} v(x) (f\sigma(y))
v(y),
\end{equation}

\noindent for all $x,\ y \in \Ga$. We denote by $Z^{1}_{f,
\sigma}(\Ga,G)$ the set of 1-cocycles associated to this action.

\smallbreak Denote by $\Emb (\Ga,G)$ the set of embeddings of
$\Ga$ in $G$ (that is, injective rational homomorphisms).

\begin{definition}\label{defequivrel}
Let $\Ga$ be an algebraic group and $\sigma_{1},\ \sigma_{2} \in
\Emb(\Ga,G)$. We say that $\sigma_{1}$ is {\it equivalent} to
$\sigma_{2}$, $\sigma_{1} \sim \sigma_{2}$, if there exist $\tau \in
\Aut(\Ga)$, $f \in \qAut(G)$ and $v \in \Map(\Ga,G)$ such that
\begin{equation}\label{eqn:spellingcocycle}
\sigma_{1}(\tau(x)) = f(\sigma_{2}(x)) v(x)\mbox{ for all }x\in
\Ga.
\end{equation}
\end{definition}

\noindent It is straightforward to see that $\sim$ is an equivalence
relation in $\Emb(\Ga,G)$. Note that the map $v$ is uniquely
determined by $f\sigma_{2}$ and $\sigma_{1}\tau $ with $v(x) =
f(\sigma_{2}(x))^{-1} \sigma_{1}(\tau(x))$ for all $x \in \Ga$; in
particular, it is a rational homomorphism. An easy computation shows
that $v \in Z^{1}_{f,\sigma_{2}}(\Ga,G)$.

\begin{obs}\label{coborde}
If the 1-cocycle $v$ is a 1-coboundary, then there exists $g \in G$
such that $v(x) = \partial (g)(x) = (g\leftharpoonup x)g^{-1} =
f(\sigma_{2}(x))^{-1} g f(\sigma_{2}(x)) g^{-1}$ for all $x\in \Ga$
and this implies that
$$\sigma_{1}(\tau(x)) = g f(\sigma_{2}(x))
g^{-1}.$$

\noindent That is, $\sigma_{1}$ can be obtained by the automorphism
$\tau,\ f$ and the conjugation by an element of $G$.
\end{obs}

\begin{definition}\label{defTfsigma}
Let $\sigma \in \Emb(\Gamma, G)$ and $f \in \qAut(G)$. Define
$\T^{f\sigma}$ to be the subgroup of $\T$ given by
$$ \T^{f\sigma} = \bigcap_{y\in \Ga} f\sigma(y)\T f\sigma(y^{-1}).$$
\end{definition}

For any $\sigma \in \Emb(\Ga, G)$ and $f \in \qAut(G)$, the
subgroup $\T^{f\sigma}$ is stable under the action of $\Ga$ on $G$
defined by $g\leftharpoonup x = f\sigma(x^{-1}) g f\sigma(x)$ for
all $x \in \Ga$ and $g \in G$. We will see in Theorem
\ref{clasesiso} that the cocycle $v$ arising in
\eqref{eqn:spellingcocycle} actually belongs to
$Z^{1}_{f,\sigma_{2}}(\Ga,\T^{f\sigma_{2}})$.

\begin{lema}\label{Z1invariante}
$Z^{1}_{f,\sigma}(\Ga,\T^{f\sigma}) = Z^{1}_{(t\cdot
f),\sigma}(\Ga,\T^{(t\cdot f)\sigma})$ for all $t\in \T$.
\end{lema}

\begin{proof} We first claim that
$\T^{(t\cdot f)\sigma} = \T^{f\sigma}$ for all $t\in \T$. Indeed,
\begin{align*}
\T^{(t\cdot f)\sigma} & = \bigcap_{y\in \Ga} (t\cdot f)\sigma(y)\T
(t\cdot f)\sigma(y^{-1}) = \bigcap_{y\in \Ga} tf\sigma(y)t^{-1}\T
tf\sigma(y^{-1})t^{-1} \\
& = \bigcap_{y\in \Ga} tf\sigma(y)\T f\sigma(y^{-1})t^{-1} =
t\T^{f\sigma} t^{-1} = \T^{f\sigma},
\end{align*}

\noindent since $\T^{f\sigma} \subset \T$ and $t\in \T$.  Suppose
that $v \in Z^{1}_{(t\cdot f),\sigma}(\Ga,\T^{(t\cdot f)\sigma})$,
then $v(x) \in \T^{(t\cdot f)\sigma}= \T^{f\sigma}$ for all $x\in
\Ga$. Now we show that $v$ is a 1-cocycle with respect to $f$ and
$\sigma$:
\begin{align*}
v(xy) & = (t\cdot f)(\sigma(y^{-1})) v(x) (t\cdot
f)(\sigma(y))v(y)\\
& = tf(\sigma(y^{-1}))t^{-1} v(x) t
f(\sigma(y))t^{-1}v(y)\\
& = tf(\sigma(y^{-1}))v(x)
f(\sigma(y))v(y)t^{-1}\\
& = f(\sigma(y^{-1}))v(x) f(\sigma(y))v(y).
\end{align*}

\noindent The last equality follows from the fact that
$f(\sigma(y^{-1}))v(x) f(\sigma(y))v(y) \in \T$ for all $x,\ y \in
\Ga$, by Definition \ref{defTfsigma}. The other inclusion follows
from similar computations.
\end{proof}

A key step in the proof of one of our main results is the next
lemma.

\begin{lema}\label{inclCenH1} Fix $\sigma \in \Emb(\Ga, G)$,
$\tau \in \Aut (\Ga)$, $f \in \qAut(G)$ and define
$$\zh_{\sigma,f,\tau} =\{\eta \in \Emb(\Ga,G)\vert\
\sigma\sim \eta \mbox{ via the triple }(\tau,f,v),\ v \in
Z^{1}_{f,\sigma}(\Ga, \T^{f\sigma}) \}.$$

\noindent Then $\T^{f\sigma}$ acts on $\zh_{\sigma, f,\tau}$ by
${t}\cdot \eta = \Int (t) \eta $, $\eta\in \zh_{\sigma, f,\tau}$
and there is a bijective map
\begin{align*}
\zh_{\sigma, f,\tau}/\T^{f\sigma} \to
H^{1}_{f,\sigma}(\Ga,\T^{f\sigma}),\qquad [\eta] \mapsto
[v_{\eta}],
\end{align*}

\noindent where $v_{\eta}(x)= f(\sigma(x))^{-1}\eta(\tau(x))$ for
all $x\in \Ga$, is the unique 1-cocycle such that $\eta(\tau(x)) =
f(\sigma(x)) v_{\eta}(x)$ for all $x\in \Ga$.
\end{lema}

\begin{proof}
Let $\eta \in \zh_{\sigma, f,\tau}$. By definition there is $v=
v_{\eta} \in Z^{1}_{f,\sigma}(\Ga, \T^{f\sigma})$ such that
$\eta(\tau(x)) = f(\sigma(x))v(x)$ for all $x\in \Ga$. Then for
all $t\in \T^{f\sigma}$ we have
\begin{align*}
({t}\cdot\eta)(\tau(x)) & = t \eta(\tau(x)) t^{-1} = t
f(\sigma(x))v(x)t^{-1}\\
& = f(\sigma(x)) f(\sigma(x))^{-1} t f(\sigma(x))v(x)t^{-1}\\
& = f(\sigma(x))[ f(\sigma(x))^{-1} t f(\sigma(x))] v(x)t^{-1}\\
& = f(\sigma(x))(t \leftharpoonup x) v(x)t^{-1} \\
& = f(\sigma(x))(t\cdot v)(x),
\end{align*}

\noindent which implies that ${t}\cdot\eta \in \zh_{\sigma,
f,\tau}$, and that $\eta = {t}\cdot\eta$ if and only if $v_{\eta}
= t\cdot v_{\eta} = v_{{t}\cdot\eta}$ for all $t\in \T^{f\sigma}$.
Hence the map $(t, \eta) \mapsto {t}\cdot\eta$ defines an action
of $\T^{f\sigma}$ on $\zh_{\sigma, f,\tau}$. Let $[\eta]$ denote
the class of $\eta$ in $\zh_{\sigma, f,\tau}/\T^{f\sigma} $. Then
$[\eta] = [\nu]$ if and only if there exists $t \in \T^{f\sigma}$
such that ${t}\cdot\eta = \nu$. Thus, the map defined by $$
\zh_{\sigma, f,\tau}/\T^{f\sigma} \to
H^{1}_{f,\sigma}(\Ga,\T^{f\sigma}),\qquad [\eta] \mapsto
[v_{\eta}],
$$

\noindent is well-defined and injective. The surjectivity follows by
definition.
\end{proof}

Recall that an algebraic group $G$ is a {\it d-group} if the
coordinate ring $\OO(G)$ has a basis consisting of characters.
Clearly any torus $T$ is a $d$-group.

\begin{lema}\label{dgrupo}\cite[Prop. 16.1 and Thm. 16.2]{Hu}
\begin{enumerate}
\item[$(a)$] If $H$ is a closed subgroup of a $d$-group $G$, then
$H$ is also a $d$-group.
\item[$(b)$] A connected $d$-group is a torus. \qed
\end{enumerate}
\end{lema}

The following lemma is also crucial for the proof of our last
theorem.
\begin{lema}\label{H1finito} If $\Ga$ is finite, then
the group $H^{1}_{f,\sigma}(\Ga,\T^{f\sigma})$ is also finite.
\end{lema}

\begin{proof}
Let $\T^{f\sigma}_{0}$ be the connected component of
$\T^{f\sigma}$ which contains the identity and let $\mathfrak{T} =
\T^{f\sigma}/ \T^{f\sigma}_{0}$; then $\vert \mathfrak{T}\vert$ is
finite. Since $\T$ is closed, it follows that $f(\sigma(x)) \T
f(\sigma(x^{-1}))$ is closed for all $x\in \Ga$, thus
$\T^{f\sigma}$ is also closed. Then by Lemma \ref{dgrupo} $(a)$,
$\T^{f\sigma}$ and consequently $\T^{f\sigma}_{0}$ are $d$-groups.
Since $\T^{f\sigma}_{0}$ is connected, it follows from Lemma
\ref{dgrupo} $(b)$ that $\T^{f\sigma}_{0}$ is a torus and
consequently $\T^{f\sigma}_{0}$ is a divisible group.

As the action of $\Ga$ on $\T^{f\sigma}$ is given by the
conjugation via $f\sigma$, every $x \in \Ga$ acts on
$\T^{f\sigma}$ by a continuous automorphism. Hence,
$\T^{f\sigma}_{0}$ is stable under the action of $\Ga$ and whence
the action of $\Ga$ on $\T^{f\sigma}$ induces an action of $\Ga$
on $\mathfrak{T}$. Thus we have an exact sequence of $\Ga$-modules
$$\xymatrix{ \T^{f\sigma}_{0} \ar@{^{(}->}[r]^{\alpha} & \T^{f\sigma}
\ar@{>>}[r]^{\beta} & \mathfrak{T},}$$

\noindent which by \cite[Prop. III.6.1]{Br} induces the long exact
sequence $$\xymatrix{ 0 \ar[r] &
H^{0}_{f,\sigma}(\Ga,\T^{f\sigma}_{0}) \ar[r]^{\alpha^{0}_{*}} &
H^{0}_{f,\sigma}(\Ga,\T^{f\sigma}) \ar[r]^{\beta^{0}_{*}} &
H^{0}_{f,\sigma}(\Ga,\mathfrak{T}) \ar[r] & {} }$$
$$\xymatrix{ {}\ar[r] &
H^{1}_{f,\sigma}(\Ga,\T^{f\sigma}_{0})
 \ar[r]^{\alpha^{1}_{*}}
& H^{1}_{f,\sigma}(\Ga,\T^{f\sigma}) \ar[r]^{\beta^{1}_{*}} &
H^{1}_{f,\sigma}(\Ga,\mathfrak{T}) \ar[r] & \cdots }$$

\noindent Since $\Ga$ and $\mathfrak{T}$ are finite groups, the
order of $H^{1}_{f,\sigma}(\Ga,\mathfrak{T})$ is finite. Hence, it
is enough to show that $\vert H^{1}_{f,\sigma}(\Ga,
\T^{f\sigma}_{0}) \vert $ is finite, because in such a case $\vert
H^{1}_{f,\sigma}(\Ga, \T^{f\sigma}) \vert$ $ = \vert \Img
\alpha^{1}_{*} \vert \vert \Img \beta^{1}_{*} \vert$ is also finite.

\smallbreak By \cite[Cor. III.10.2]{Br}, we know that
$H^{n}_{f,\sigma}(\Ga,\T^{f\sigma}_{0})$ is annihilated by $m= \vert
\Ga \vert $  for all $n> 0$. Then, for all $\alpha \in
Z^{1}_{f,\sigma}(\Ga, \T^{f\sigma}_{0})$ there exists $t \in
\T^{f\sigma}_{0}$ such that $\alpha^{m} = \partial (t)$. Since
$\T^{f\sigma}_{0}$ is divisible, there exists $s \in
\T^{f\sigma}_{0}$ such that $t = s^{m}$. Let $\beta =
\partial(s^{-1}) \alpha$, then $\beta^{m} = 1$ and therefore
$\beta \in Z^{1}_{f,\sigma}(\Ga, D_{m})$, where $D_{m} = \{ t \in
\T^{f\sigma}_{0}\vert\ t^{m} =1\}$. Moreover, $[\alpha] = [\beta]$
and the inclusion of the set of 1-cocycles defines a surjective map
$H^{1}_{f,\sigma}(\Ga, D_{m})
 \to H^{1}_{f,\sigma}(\Ga,\T^{f\sigma}_{0})$. As $\T^{f\sigma}_{0}$
is a torus,  $D_{m}$ is a finite group and consequently $\vert
H^{1}_{f,\sigma}(\Ga, D_{m})\vert $ is finite, which implies that
$\vert H^{1}_{f,\sigma}(\Ga,\T^{f\sigma}_{0})\vert$ is also
finite.
\end{proof}

\subsection{Isomorphisms}
In this subsection we study the isomorphisms between the Hopf
algebras $A_{\sigma}$ for a fixed algebraic group $\Ga$.

\smallbreak By Corollary \ref{propL}, the $\qe^{*}$-extension $\Oe$
of $\OO(G)$ satisfies the property $(L)$. Moreover, by \cite[Prop.
A.3]{AS-app}, the Frobenius-Lusztig kernels $\qe$ are simple Hopf
algebras if the order $\ell$ of the root of unity $\epsilon$ and
$\det DC$ are relatively prime, where $DC$ is the ``symmetrized''
Cartan matrix of $\lieg$. If $\lieg$ is not of type $A_{n}$, this is
always the case. In this situation $\Z(\qe^{*}) = \CC$, since
otherwise $\qe^{*}$ would contain a central proper Hopf subalgebra
${\bf v}$ and thus $\qe^{*}$ is an extension of ${\bf v}$ by
$\qe^{*} / {\bf v}^{+}\qe^{*}$, implying by duality that $\qe$
contains a proper normal Hopf subalgebra dual to $\qe^{*} / {\bf
v}^{+}\qe^{*}$.

\smallbreak Hence, $\qe^*$ satisfies the property $(Z)$ if  $\ell$
and $\det DC$ are relatively prime. For this reason, we assume
from now on that $$ \ell \text{ and } \det DC \text{ are
relatively prime.}
$$

\begin{obs} The Hopf algebras $A_{\sigma} = \Oe /(\JJ)$ not only
depend on $\sigma$ but also on $\lieg$ and the root of unity
$\epsilon$. Nevertheless, the last data are invariant with respect
to the first construction: if $A_{\sigma, \epsilon, \lieg} \simeq
A_{\sigma', \epsilon', \lieg'}$ then $\epsilon = \epsilon'$, $\lieg
\simeq \lieg'$ and $\Ga \simeq \Ga'$. Indeed, let $\varphi:
A_{\sigma, \epsilon, \lieg} \to A_{\sigma', \epsilon', \lieg'}$
denote the isomorphism. Since $\Z(\qe^{*}) = \Z(\qep^{*}) = \CC$, by
Corollary \ref{sincentro} we have that $\Z(A_{\sigma, \epsilon,
\lieg}) = \Oc(\Ga)$ and $\Z(A_{\sigma', \epsilon', \lieg'}) =
\Oc(\Ga')$. Thus by Proposition \ref{isoinducidos}, $\varphi$
induces the isomorphisms $\underline{\varphi}: \Oc(\Ga) \to
\Oc(\Ga')$ and $\overline{\varphi}: \qe^{*} \to \qep^{*}$. In
particular, we have that $\Ga \simeq \Ga'$ and $\qe \simeq \qep$.
Since the Frobenius-Lusztig kernels are quantum linear spaces of
finite Cartan type, from \cite[Thm. 7.2]{AS-05} it follows that
$\epsilon = \epsilon'$ and $\lieg \simeq \lieg'$.
\end{obs}

Fix the simple Lie algebra $\lieg$, the root of unity $\epsilon$ and
an algebraic group $\Ga$. Let $\sigma_{1},\ \sigma_{2} \in
\Emb(\Ga,G)$ and denote $A_{i}= A_{\sigma_{i}}$. Recall that by
Definition \ref{defequivrel}, $\sigma_{1}\sim \sigma_{2}$ if  there
is a triple $(\tau, f, v)$ such that
\begin{enumerate}
\item[$(i)$] $\tau \in \Aut(\Ga)$,
\item[$(ii)$] $f \in \qAut(G)$ and
\item[$(iii)$] $v \in Z^{1}_{f,\sigma_{2}}(\Ga,G)$
with $v(1) = 1$.
\end{enumerate}
Now we are able to apply Theorem \ref{isoytriple}, which in this
case gives:

\begin{theorem}\label{clasesiso} $A_{1}$ and $A_{2}$ are
isomorphic if and only if $\sigma_{1}\sim \sigma_{2}$ via a triple $(\tau,
f, v)$ with $v \in Z^{1}_{f,\sigma_{2}}(\Ga,\T^{f\sigma_{2}})$.
\end{theorem}

\begin{proof} By Theorem \ref{isoytriple},
we know that the Hopf algebras $A_{1}$ and $A_{2}$ are isomorphic
if and only if there is a triple $(\underline{\omega},g, u)$ such
that $\underline{\omega} \in \Aut(\Oc(\Ga))$, $g \in
\qAut(\OO(G))$ and $u: \Oe \to \Oc(\Ga)$ is an algebra map
satisfying \eqref{conmutconu} and \eqref{comdeu}. The transposes
of these maps induce maps $^{t}\underline{\omega} = \tau \in
\Aut(\Ga)$, $^{t}g = f\in \qAut(G)$ and $^{t}u = \mu \in \Map(\Ga,
\Ge)$, where $\Ge = \Alg(\Oe,\CC)$ and $\mu(1)= 1$.

Let $v =\ ^{t}\iota\mu$, then $v$ is a 1-cocycle which satisfies
$(iii)$: by Lemma \ref{lema:toro}, the image of $\ ^{t}\iota: \Ge
\to G$ is the maximal torus $\T$, thus  $v =\ ^{t}\iota \mu$ is a
map $v:\Ga \to \T \subset G$, which by \eqref{conmutconu} and
\eqref{comdeu} satisfies that

\begin{align*}
\sigma_{1}(\tau(x)) & = f(\sigma_{2}(x)) v(x), \mbox{ and }
\end{align*}

\begin{align*}
v(x y) & =  f(\sigma_{2}(y))^{-1}v(x)\sigma_{1}(\tau(y))\\
& = f(\sigma_{2}(y))^{-1}v(x)f(\sigma_{2}(y)) v(y)\\
& = (v(x)\leftharpoonup y)v(y).
\end{align*}

\noindent for all $g,\ h \in \Ga$. Moreover, from the equalities
above it follows that
$$v(xy) v(y)^{-1} = f(\sigma_{2}(y))^{-1}v(x)f(\sigma_{2}(y))
\mbox{ for all }x,\ y \in \Ga,$$

\noindent which implies that
$f(\sigma_{2}(y))^{-1}v(x)f(\sigma_{2}(y)) \in \T$ for all $x,\ y
\in \Ga$, that is, $v(x) \in \T^{f\sigma_{2}}$ for all $x\in \Ga$,
and thus $v \in Z^{1}_{f,\sigma_{2}}(\Ga, \T^{f\sigma_{2}})$. By
definition, both equalities above hold if and only if $\sigma_{1}
\sim \sigma_{2}$ via the maps $\tau \in \Aut(\Ga)$, $f\in
\qAut(G)$ and $v \in Z^{1}_{f,\sigma_{2}}(\Ga, \T^{f\sigma_{2}})$.
\end{proof}

Observe that $\Aut(G)$ acts on $\Emb(\Ga,G)$ by $f\circ \sigma$
for all $f \in \Aut(G)$ and $\sigma \in \Emb(\Ga,G)$. In
particular, $G$ and $\T$ act on $\Emb(\Ga,G)$ by $\Int(G)$ and
$\Int(\T)$ respectively. Denote $\Int(g)\circ \sigma ={g}\cdot
\sigma$, for all $\sigma \in \Emb(\Ga,G)$, $g\in G$ and $
G\cdot\sigma$ the orbit of $\sigma$ in $\Emb(\Ga,G)$ under the
action of $G$. Clearly,  if $C = C_{G}(\sigma(\Ga))$ is the
centralizer of $\sigma(\Ga)$ in $G$, then $G\cdot\sigma\simeq G/C$
and the set of $\T$-orbits in $G\cdot\sigma$ is $\T\backslash
G/C$.

\bigbreak
\emph{Assume from now on that $\Ga$ is finite.}

\begin{lema}\label{dimgmayrg}
If $\sigma(\Ga)$ is not central in $G$, then $\T\backslash G/C$ is
infinite.
\end{lema}

\begin{proof}
To prove that $\T\backslash G/C$ is infinite it is enough to show
that $\dim G > \dim \T + \dim C =  \rk G + \dim C$, since if $m =
\# \T\backslash G/C$ were finite, then $G = \bigcup_{i=1}^{m} \T
g_{i} C$ and this would imply that $\dim G \leq \dim \T + \dim C$.
Pick $g\in \sigma(\Ga)$ non-central; then the centralizer of $g$
is a proper reductive subgroup of $G$ by \cite{R}, and $C$ is
contained in a maximal reductive subgroup $M$ of $G$. Let $\lieg$
and $\liem$ be the Lie algebras of $G$ and $M$ respectively. As
the maximal subalgebras of the simple Lie algebras are classified,
by inspection in \cite{D1,D2} one can see that $\dim \lieg
> \rk \lieg + \dim \liem $ for any maximal reductive $\liem$. See
\cite[App.]{G} for details.
\end{proof}

\begin{lema}\label{fliadualnopunt}
If $\sigma(\Ga)$ is not central in $G$, then there exists an
infinite family $\{g_{i}\}_{i\in I}$ of elements of $G$ such that
$(g_{i}\cdot \sigma) (\Ga) \nsubseteq \T$.
\end{lema}

\begin{proof}
Let $\{g_{j}\}_{j\in J}$ be a set of representatives of
$\T\backslash G/C$. Then, $J$ is infinite by Lemma \ref{dimgmayrg}.
We will prove that there can only be finitely many $g_{j}$ such that
$(g_{j}\cdot \sigma) (\Ga) \subseteq \T$.

\smallbreak Suppose that there exists $g_{i} \in G$ such that
$(g_{i}\cdot \sigma)(\Ga) \subseteq \T$. Without loss of
generality we can assume that $\Ga \subseteq \T$. Consider the
sets
\begin{align*}
\L  = \{h \in G\vert\ h\Ga h^{-1} \subseteq \T\} \qquad \mbox{ and }
\qquad \Lambda = \sum_{h \in \L} h\Ga h^{-1} .
\end{align*}

\noindent Then clearly $\T \subseteq \L$ and $\Lambda \subset \T$.
We show that $\T \backslash \L / C $ is finite. Let $N = \vert \Ga
\vert$, then $N \Lambda = 0$, which implies that $\Lambda
\subseteq (\GG_{N})^{n}$, where $\GG_{N}$ is the  group of $N$-th
roots of unity. Thus $\Lambda$ is a finite subgroup of $\T$. In
particular, it contains only finitely many subgroups which are
pairwise distinct and isomorphic to $\Ga$. Let $\Ga_{1},\ldots
,\Ga_{s}$ be these subgroups and $h_{i} \in \L$ such that $h_{i}
\Ga h_{i}^{-1} = \Ga_{i}$. If $h \in \L$, then $h\Ga h^{-1} =
h_{i} \Ga h_{i}^{-1}$ for some $1\leq i\leq s$. Hence $h_{i}^{-1}h
\Ga h^{-1}h_{i} = \Ga$ and $h_{i}^{-1}h \in N_{G}(\Ga)$, the
normalizer of $\Ga$ in $G$. Thus $\L = \coprod_{i=1}^{s}
h_{i}N_{G}(\Ga)$.

\smallbreak On the other hand, there is a homomorphism $N_{G}(\Ga)
\to \Aut(\Ga)$ which factorizes through the monomorphism
$N_{G}(\Ga) / C \hookrightarrow \Aut(\Ga)$. Since $\Ga$ is finite,
the order of $N_{G}(\Ga) / C$ is finite and consequently $\L / C$
is finite. Since by assumption $\Ga \subseteq \T$, there exists an
injection $\T \backslash \L / C \hookrightarrow \L / C$, which
implies that $\{g_{j}\}_{j\in J}$ contains only finitely many
elements such that $g_{j}\Ga g_{j}^{-1} \subseteq \T$.
\end{proof}

Now we are able to prove our last theorem.

\begin{theorem}\label{fliainfinita}
Let $\sigma \in \Emb(\Ga,G)$ such that $\sigma(\Ga)$ is not central
in $G$. Then there exists an infinite family $\{\sigma_{j}\}_{j\in
J} \subset \Emb(\Ga, G)$ such that the Hopf algebras
$\{A_{\sigma_{j}}\}_{j\in J}$ of dimension $\vert \Ga \vert
\ell^{\dim \lieg}$ are pairwise non-isomorphic, non-semisimple,
non-pointed and their duals are also non-pointed.
\end{theorem}

\begin{proof}
Since $\sigma(\Ga)$ is not central in $G$, from Lemma
\ref{dimgmayrg} it follows that there exists an infinite set
$\{g_{i}\}_{i\in I}$ of elements in $G$ such that
$({tg_{i}})\cdot\sigma \neq\ {g_{j}}\cdot\sigma$ for all $i\neq j$
and $t\in \T$. Denote $\sigma_{i} = g_{i} \cdot\sigma$ for all $i\in
I$. By Definition \ref{defequivrel} we know that
$$ \Emb(\Ga, G ) = \coprod_{\eta \in E} \zh_{\eta},$$

\noindent where $\zh_{\eta} = \{\mu \in \Emb(\Ga,G)\vert\ \eta\sim
\mu \}$ and $E$ is a set of representatives of $\Emb(\Ga,G)$ under
the equivalence relation $\sim$. We prove now that there can not
be infinitely many embeddings $\sigma_{i}$ included in only one
$\zh_{\eta}$.

\smallbreak Suppose on the contrary, that there exists $\eta \in
\Emb(\Ga,G)$ such that $\zh_{\eta}$ contains infinitely many
$\sigma_{i}$. By definition and Theorem \ref{clasesiso}, we know
that $$ \zh_{\eta} = \bigcup_{\tau\in \Aut (\Ga),\, f \in \qAut(G)}
\zh_{\eta, f,\tau},$$

\noindent where $\zh_{\eta, f,\tau} =\{\mu \in \Emb(\Ga,G)\vert\
\eta\sim \mu \mbox{ via  }(\tau,f,v)$, with $v \in
Z^{1}_{f,\eta}(\Ga, \T^{f\eta})\}$.

\smallbreak By Lemma \ref{TactuaenqAut}, $\T$ and $\T^{f\eta}$ act
on $\qAut(G)$ by $t\cdot f (x) = tf(x)t^{-1}$ for all $f \in
\qAut(G)$, $t \in \T$. We claim that if $t \in \T^{f\eta}$, then
$\zh_{\eta, t\cdot f, \tau} = \zh_{\eta, f, \tau}$. Indeed, $\mu
\in \zh_{\eta, t\cdot f, \tau}$ if and only if there exists $v \in
Z^{1}_{(t\cdot f),\eta}(\Ga, \T^{(t\cdot f)\eta})$ such that
$\mu(\tau(x)) = (t\cdot f)(\eta(x)) v(x)$, for all $x\in \Ga$. But
in such a case
\begin{align*}
\mu(\tau(x))  & = (t\cdot f)(\eta(x)) v(x) = tf(\eta(x))
t^{-1}v(x) = tf(\eta(x)) v(x) t^{-1}\\
& = f(\eta(x)) [f(\eta(x)) ^{-1}tf(\eta(x))] v(x) t^{-1} =
f(\eta(x)) (t\cdot v)(x),
\end{align*}

\noindent which implies that $\mu \in \zh_{\eta, f, \tau}$, since
by Lemma \ref{Z1invariante}, $Z^{1}_{(t\cdot f),\eta}(\Ga,
\T^{(t\cdot f)\eta}) = Z^{1}_{f ,\eta}(\Ga, \T^{f\eta})$ for all
$t\in \T$, and by Lemma \ref{inclCenH1}, $t\cdot v \in Z^{1}_{f
,\eta}(\Ga, \T^{f\eta})$. Thus, $\zh_{\eta, t\cdot f, \tau}
\subseteq \zh_{\eta, f, \tau}$ and similarly for the other
inclusion. Hence we can write
$$ \zh_{\eta} = \coprod_{\tau\in \Aut(\Ga)}
\coprod_{f \in \J} \coprod_{t \in \t} \zh_{\eta, t\cdot f,\tau},$$

\noindent where $\J$ is a set of representatives of $\qAut(G)/\T$
and $\t$ is a set of representatives of $\T / \T^{f\eta}$. Since
$\Aut(\Ga)$ is finite and by Lemma \ref{TactuaenqAut} (c), $\J$ is
also finite, if $\zh_{\eta}$ contains infinitely many $\sigma_{i}$
then there exist $\tau \in \Aut(\Ga)$ and $f \in \J$ such that
$\coprod_{t \in \t} \zh_{\eta, t\cdot f,\tau}$ contains infinite
many $\sigma_{i}$. If $\sigma_{i} \in \zh_{\eta, t\cdot f,\tau}$ for
some $t \in \t$, then ${t^{-1}}\cdot\sigma_{i} \in \zh_{\eta,
f,\tau}$. By Lemmas \ref{inclCenH1} and \ref{H1finito}, we know that
the set $\zh_{\eta, f,\tau} / \T^{f\eta}$ is finite, hence there
must exists $\sigma_{j}$, $i\neq j$ and $s \in \T$ such that
$[t^{-1}\cdot\sigma_{i}] = [s^{-1}\cdot\sigma_{j}]$. But this
contradicts our assumption on the family $\{\sigma_{i}\}_{i\in I}$
since in such a case, there would exist $r \in \T^{f\eta}$ such that
${t^{-1}}\cdot\sigma_{i} =\ {r}\cdot({s^{-1}}\cdot\sigma_{j}) =\
({rs^{-1}})\cdot\sigma_{j}$, that is $\sigma_{i} =\
({trs^{-1}})\cdot\sigma_{j}$ with $trs^{-1} \in \T$.

\smallbreak In conclusion, there can be only finitely many
elements of the set $\{\sigma_{i}\}_{i\in I}$ in each
$\zh_{\eta}$, for any $\eta \in E$. Thus, by Corollary
\ref{cor:pulenta}, Theorem \ref{clasesiso} and Lemma
\ref{fliadualnopunt}, there exists an infinite subset $J\subset I$
such that the Hopf algebras $\{A_{\sigma_{j}}\}_{j\in J}$ are
pairwise non-isomorphic, non-semisimple, non-pointed and with
non-pointed duals.
\end{proof}

\begin{obs}
$(a)$ The Hopf algebras in the preceding theorem are all
quasi-isomorphic. This follows from \cite{Mk}, for they are
constructed via a pushout. Indeed, if $\Ga$ and $\Ga'$
are conjugate in $G$, then $A_{\Ga}$ and $A_{\Ga'}$ are
cocycle deformations of each other.

\smallbreak $(b)$ The examples given by M\" uller in \cite[Thm.
5.13]{Mu2} fit into the situation above for $\lieg =
\mathfrak{sl}_{2}$ and $\Ga = \GG_{\ell}$.

\smallbreak $(c)$ The finite subgroups of a simple affine
algebraic group $G$ are not known. However, some positive results
for subgroups of prime order were proved by several authors. See
for example \cite{CG}, \cite{PW}, \cite{S}.

\end{obs}


\subsection*{Acknowledgments} The authors thank S. Natale
for helpful indications on Lemma \ref{H1finito} and J. Vargas for
discussions on reductive subgroups of a simple group. G. A. G.
thanks H.-J. Schneider for hospitality and fruitful conversations
during his stay in Munich. They also thank the referees for a very
careful reading of the paper.


\bibliographystyle{amsbeta}

\end{document}